\documentclass[a4paper,11pt]{article}
\usepackage{amsmath,amsthm,amssymb,enumitem,xcolor}

\usepackage[hang]{footmisc}
\usepackage[nosort,nocompress,noadjust]{cite}

\usepackage[bookmarks=false,hyperfootnotes=false,colorlinks,
    linkcolor={red!60!black},
    citecolor={blue!50!black},
    urlcolor={blue!80!black}]{hyperref}

\renewcommand{\eqref}[1]{\hyperref[#1]{(\ref{#1})}}

\newlist{enumlist}{enumerate}{2}
\setlist[enumlist,1]{labelindent=0cm,label=\arabic*.,ref=\arabic*,labelwidth=2.5ex,labelsep=0.5ex,leftmargin=3ex,align=left,topsep=0.5ex,itemsep=1ex,parsep=1ex}
\setlist[enumlist,2]{labelindent=0cm,label=\theenumlisti.\arabic*.,ref=\arabic*,labelwidth=5ex,labelsep=0.5ex,leftmargin=5.5ex,align=left,topsep=0.5ex,itemsep=1ex,parsep=1ex}

\newlist{itemlist}{itemize}{1}
\setlist[itemlist]{labelindent=0cm,label=$\bullet$,labelwidth=2.5ex,labelsep=0.5ex,leftmargin=3ex,align=left,topsep=0.5ex,itemsep=1ex,parsep=1ex}

\numberwithin{equation}{section}

{\theoremstyle{definition}\newtheorem{definition}{Definition}[section]
\newtheorem{remark}[definition]{Remark}
\newtheorem{example}[definition]{Example}}

\newtheorem{proposition}[definition]{Proposition}
\newtheorem{lemma}[definition]{Lemma}
\newtheorem{theorem}[definition]{Theorem}
\newtheorem{corollary}[definition]{Corollary}
\newtheorem{letterthm}{Theorem}
\newtheorem{lettercor}[letterthm]{Corollary}

\newcommand{\om}{\omega}
\newcommand{\si}{\sigma}
\newcommand{\actson}{\curvearrowright}
\newcommand{\R}{\mathbb{R}}
\newcommand{\N}{\mathbb{N}}

\newcommand{\cF}{\mathcal{F}}
\newcommand{\cP}{\mathcal{P}}
\newcommand{\eps}{\varepsilon}

\newcommand{\cU}{\mathcal{U}}

\newcommand{\vphi}{\varphi}
\newcommand{\mutil}{\widetilde{\mu}}
\newcommand{\gammatil}{\widetilde{\gamma}}
\newcommand{\cR}{\mathcal{R}}
\newcommand{\cG}{\mathcal{G}}

\newcommand{\Z}{\mathbb{Z}}

\newcommand{\al}{\alpha}
\newcommand{\be}{\beta}

\newcommand{\cS}{\mathcal{S}}
\newcommand{\C}{\mathbb{C}}
\newcommand{\cN}{\mathcal{N}}

\newcommand{\Aut}{\operatorname{Aut}}
\newcommand{\Xtil}{\widetilde{X}}
\newcommand{\Om}{\Omega}

\newcommand{\Omtil}{\widetilde{\Omega}}
\newcommand{\id}{\mathord{\text{\rm id}}}
\newcommand{\Out}{\operatorname{Out}}
\newcommand{\Ker}{\operatorname{Ker}}
\newcommand{\alh}{\widehat{\alpha}}
\newcommand{\Yh}{\widehat{Y}}
\newcommand{\etah}{\widehat{\eta}}
\newcommand{\SL}{\operatorname{SL}}
\newcommand{\GL}{\operatorname{GL}}
\newcommand{\cA}{\mathcal{A}}
\newcommand{\cV}{\mathcal{V}}
\newcommand{\cH}{\mathcal{H}}
\newcommand{\Q}{\mathbb{Q}}
\newcommand{\yh}{\widehat{y}}
\newcommand{\beh}{\widehat{\beta}}
\newcommand{\sitil}{\widetilde{\sigma}}
\newcommand{\Deltatil}{\widetilde{\Delta}}
\newcommand{\cO}{\mathcal{O}}
\newcommand{\Cl}{\operatorname{Cl}}
\newcommand{\Gtil}{\widetilde{G}}
\newcommand{\omtil}{\widetilde{\omega}}
\newcommand{\Gammatil}{\widetilde{\Gamma}}

\newcommand{\cB}{\mathcal{B}}
\newcommand{\etil}{\widetilde{e}}
\newcommand{\Ctil}{\widetilde{C}}
\newcommand{\Ttil}{\widetilde{T}}

\begin{document}

\begin{center}
{\boldmath\LARGE\bf Orbit equivalence superrigidity for type III$_0$ actions}

\bigskip

{\sc by Stefaan Vaes\footnote{\noindent KU~Leuven, Department of Mathematics, Leuven (Belgium).\\ E-mails: stefaan.vaes@kuleuven.be and bram.verjans@kuleuven.be.}\textsuperscript{,}\footnote{S.V.\ is supported by FWO research project G090420N of the Research Foundation Flanders and by long term structural funding~-- Methusalem grant of the Flemish Government.} and Bram Verjans\textsuperscript{1,}\footnote{B.V.\ is holder of PhD grant 1139721N fundamental research of the Research Foundation~-- Flanders.}}
%
%
\end{center}

\begin{abstract}
\noindent We prove the first orbit equivalence superrigidity results for actions of type III$_\lambda$ when $\lambda \neq 1$. These actions arise as skew products of actions of dense subgroups of $\SL(n,\R)$ on the sphere $S^{n-1}$ and they can have any prescribed associated flow.
\end{abstract}

\section{Introduction and main results}

An essentially free, ergodic, nonsingular action $G \actson (X,\mu)$ of a countable group $G$ on a standard probability space is said to be orbit equivalence (OE) superrigid if the group $G$ and its action on $(X,\mu)$ can be entirely retrieved from the orbit equivalence relation $\cR(G \actson X) = \{(x,g \cdot x) \mid x \in X, g \in G\}$. Especially in the case where $\mu$ is a $G$-invariant probability measure, several OE-superrigidity theorems were proven in the context of Popa's deformation/rigidity theory, see e.g.\ \cite{Pop05,Pop06,Ioa08,DIP19}.

Beyond the probability measure preserving setting, OE-superrigidity results are more scarce, see e.g.\ \cite{PV08,Ioa14,DV21}. In all these cases, the group action $G \actson (X,\mu)$ has one of the following Krieger types: II$_1$, II$_\infty$ or III$_1$ (see Section \ref{sec.prelim} for definitions). There is a conceptual reason why it is harder to prove OE-superrigidity for actions of type III$_\lambda$, $\lambda \in [0,1)$. One typically proves OE-superrigidity for $G \actson (X,\mu)$ by showing that every measurable $1$-cocycle $\om : G \times X \to \Lambda$ with values in an arbitrary countable group $\Lambda$ is cohomologous to a group homomorphism $\delta : G \to \Lambda$. When the measure $\mu$ is not $G$-invariant, the logarithm of the Radon-Nikodym derivative $d(g^{-1}\cdot \mu)/d\mu$ provides a $1$-cocycle $\om : G \times X \to \R$. In type III$_\lambda$ with $\lambda \in [0,1)$, this $1$-cocycle is ``essentially'' similar to a $1$-cocycle with values in a countable group. Therefore, cocycle superrigidity tends to fail.

In this paper, we obtain the first OE-superrigidity results in type III$_\lambda$ when $\lambda \neq 1$. In Theorem~\ref{thm.main-no-go-theorem}, we prove OE-superrigidity for the affine action of dense subgroups $G < \SL(n,\R) \ltimes \R^n$ on $X = \R^n$. These actions can be of type III$_\lambda$ for any $\lambda \in (0,1]$. In this result, OE-superrigidity holds in its strongest possible form: for every essentially free, ergodic, nonsingular action $\Lambda \actson (Z,\zeta)$ that is stably orbit equivalent with $G \actson (X,\mu)$, there exists an injective group homomorphism $\delta : G \to \Lambda$ such that $\Lambda \actson (Z,\zeta)$ is isomorphic to the induction of $G \actson (X,\mu)$ along $\delta$.

In Theorem \ref{thm.main-no-go-theorem}, we also prove that this strongest form of OE-superrigidity can basically never hold for actions of type III$_0$. In type III$_0$, it is necessary to further reduce the other action $\Lambda \actson (Z,\zeta)$~: after writing $\Lambda \actson Z$ as an induction of $\Lambda_0 \actson Z_0$, we need to take a quotient $\Lambda_0 / \Sigma \actson Z_0 / \Sigma$ by a normal subgroup $\Sigma$ whose action on $Z_0$ admits a fundamental domain, before arriving at an action that is conjugate with $G \actson (X,\mu)$.

This leads to a second, slightly weaker version of OE-superrigidity that we denote, without too much inspiration, as OE superrigidity (v2). In Theorem \ref{thm.main-OE-superrigid} and Corollary \ref{cor.main-example}, we then prove that natural skew product actions of dense subgroups of $\SL(n,\R)$ are OE-superrigid (v2), of type III$_0$, with any prescribed associated flow.

Before stating our main results, we make this terminology more precise. First note that the concepts of (stable) orbit equivalence, induced actions, conjugate actions, etc., are recalled in Section \ref{sec.prelim}.


Beyond the probability measure preserving (pmp) setting, one cannot distinguish between orbit equivalence and stable orbit equivalence. Therefore, induced actions will appear in any OE-superrigidity statement. We thus formally define the following property for a free, ergodic, nonsingular action $G \actson (X,\mu)$.

\begin{itemlist}
\item {\bf OE superrigidity (v1)} of $G \actson (X,\mu)$~: any free, ergodic, nonsingular action that is stably orbit equivalent with $G \actson (X,\mu)$ is conjugate to an induction of $G \actson (X,\mu)$.
\end{itemlist}

As we prove in Theorem \ref{thm.main-no-go-theorem}, this (v1) of OE superrigidity can basically never hold for actions of type III$_0$, but does hold for several actions of type III$_\lambda$ with $\lambda \in (0,1]$. For actions of type III$_0$, unavoidably the following extra freedom is needed, accommodating for the canonical stable orbit equivalences that come with induction and with quotients by normal subgroups whose action admits a fundamental domain.

\begin{itemlist}
\item {\bf OE superrigidity (v2)} of $G \actson (X,\mu)$~: if a free, ergodic, nonsingular action $\Lambda \actson (Z,\zeta)$ is stably orbit equivalent with $G \actson (X,\mu)$, there exist subgroups $\Sigma < \Lambda_0 < \Lambda$ and a nonnegligible $Z_0 \subset Z$ such that $\Lambda \actson Z$ is induced from $\Lambda_0 \actson Z_0$, $\Sigma \lhd \Lambda_0$ is normal, the action $\Sigma \actson Z_0$ admits a fundamental domain and $G \actson (X,\mu)$ is conjugate with $\Lambda_0/\Sigma \actson Z_0/\Sigma$.
\end{itemlist}

To obtain technically less involved statements, one may restrict to \emph{simple} actions: in Definition \ref{def.simple-action}, we say that a free, ergodic, nonsingular action $G \actson (X,\mu)$ is simple if the action is not induced and if $G$ has no nontrivial normal subgroups whose action on $(X,\mu)$ admits a fundamental domain. Then both versions of OE superrigidity for a simple action $G \actson (X,\mu)$ immediately imply that any stably orbit equivalent simple action must be conjugate to $G \actson X$, bringing us back to a statement that looks similar to the probability measure preserving setting.

%
%

We use the following skew product construction to obtain OE-superrigid actions of type III$_0$. Given any nonsingular ergodic action $G \actson (X,\mu)$ of type III$_1$, with logarithm of the Radon-Nikodym cocycle denoted as $\om : G \times X \to \R$, and given any ergodic flow $\R \actson^\al (Y,\eta)$, we consider
\begin{equation}\label{eq.our-actions}
G \actson (X \times Y, \mu \times \eta) : g \cdot (x,y) = (g \cdot x, \al_{\om(g,x)}(y)) \; .
\end{equation}
We prove in Proposition \ref{prop.prescribed-flow} that this action is ergodic and that its associated flow is given by the adjoint flow $\alh$, a new concept that we introduce in Definition \ref{def.adjoint-flow}. Since the adjoint operation is involutive, meaning that the adjoint of $\alh$ is isomorphic with $\al$, the skew product construction \eqref{eq.our-actions} provides a streamlined way of defining group actions with a prescribed associated flow.

The main result of this paper is the following OE superrigidity theorem for actions of type III$_0$.


\begin{letterthm}\label{thm.main-OE-superrigid}
Let $G \actson (X,\mu)$ be a free, ergodic, nonsingular action of type III$_1$. Assume that $G$ is finitely generated and has trivial center. Assume that the Maharam extension of $G \actson (X,\mu)$ is simple and cocycle superrigid with countable targets. Denote by $\om : G \times X \to \R$ the logarithm of the Radon-Nikodym cocycle.

For any ergodic flow $\R \actson^\al (Y,\eta)$, the action $G \actson (X \times Y,\mu \times \eta)$ defined in \eqref{eq.our-actions} is OE superrigid (v2) and has associated flow $\alh$.
\end{letterthm}

We provide a more precise version of Theorem \ref{thm.main-OE-superrigid} as Corollary \ref{cor.OE-superrigid-prescribed-flow} below. In this more precise version, the possible group actions $\Lambda_0 \actson Z_0$ with normal subgroup $\Sigma \lhd \Lambda_0$ that appear in the definition of OE superrigidity (v2) are explicitly described.

As we explain in Example \ref{ex.the-ex-from-PV} and Theorem \ref{thm.rings-cocycle-superrigid}, there are many concrete type III$_1$ actions $G \actson (X,\mu)$ satisfying the assumptions of Theorem \ref{thm.main-OE-superrigid}. In particular, we obtain the following result.


Recall that given a commutative ring $\cA$ and an integer $n \geq 2$, the group $E(n,\cA)$ is the subgroup of $\SL(n,\cA)$ generated by the elementary matrices having $1$'s on the diagonal and an element of $\cA$ as an off diagonal entry. For several rings, including $\Z[\cS^{-1}]$ where $\cS$ is a finite set of prime numbers and the ring of integers $\cO_K$ of an algebraic number field, we have that $E(n,\cA) = \SL(n,\cA)$ (see Example \ref{ex.allerlei-voorbeelden} for references and more examples).

\begin{lettercor}\label{cor.main-example}
Let $n \geq 3$ be an odd integer and let $\cA \subset \R$ be a subring containing an algebraic number that does not belong to $\Z$. Assume that $\cA$ is finitely generated as a ring. For every ergodic flow $\R \actson^\al (Y,\eta)$, consider the action
\begin{multline*}
\quad\quad\be(n,\cA,\al) : E(n,\cA) \actson (\R^n \times Y)/\R : A \cdot \overline{(x,y)} = \overline{(A x,y)} \\
\text{where}\quad \R \actson \R^n \times Y : t \cdot (x,y) = (e^{t/n} x,t \cdot y) \; .\quad\quad
\end{multline*}
\begin{enumlist}
\item The actions $\be(n,\cA,\al)$ are essentially free, ergodic, simple and OE superrigid (v2), with associated flow $\alh$.
\item The actions $\be(n,\cA,\al)$ and $\be(n',\cA',\al')$ are stably orbit equivalent if and only if $n = n'$, $\cA = \cA'$ and $\al$ is isomorphic with $\al'$.
\end{enumlist}
\end{lettercor}

We prove Corollary \ref{cor.main-example} as Corollary \ref{cor.main-example-with-out} below, in which we also describe the outer automorphism group $\Out(\cR(n,\cA,\al))$ of the orbit equivalence relations $\cR(n,\cA,\al)$ of the actions $\be(n,\cA,\al)$ appearing in Corollary \ref{cor.main-example}. In Remark \ref{rem.a-lot}, we also show that the family of group actions $\be(n,\cA,\al)$ in Corollary \ref{cor.main-example} is large and complex in a descriptive set theoretic sense of the word.

As mentioned above, OE superrigidity (v1) is impossible for actions of type III$_0$, but does happen for actions of type III$_\lambda$ when $\lambda \in (0,1]$. The precise result goes as follows and provides the first examples of OE superrigidity (v1) for actions of type III$_\lambda$ with $\lambda \in (0,1)$. Examples of type III$_1$ were given before, see e.g.\ \cite[Theorem 5.8]{PV08} and \cite[Proposition 3.3]{DV21}, and see Corollary \ref{cor.an-example-of-OE-superrigid-v1}.

\begin{letterthm}\label{thm.main-no-go-theorem}
\begin{enumlist}
%
\item Let $n \geq 3$ be an integer, $p$ a prime number and $0 < \lambda < 1$. Consider the ring $\cA = \Z[\lambda,\lambda^{-1},p^{-1}]$. Define the subgroup $\Gamma < \GL(n,\cA)$ of matrices $A$ with $\det A \in \lambda^\Z$.

The action of $\Gamma \ltimes \cA^n$ on $\R^n$ by $(A,a) \cdot x = A(a+x)$ is essentially free, ergodic, simple, of type III$_\lambda$. It is OE superrigid (v1).

\item Let $G \actson (X,\mu)$ be any essentially free, ergodic, simple, type III$_0$ action. Then $G \actson (X,\mu)$ is not OE superrigid (v1).
\end{enumlist}
\end{letterthm}

We prove Theorem \ref{thm.main-no-go-theorem} at the end of Section \ref{sec.versions-OE-superrigidity}.

For every free, ergodic, nonsingular action $G \actson (X,\mu)$ of a countable group $G$, the crossed product $M = L^\infty(X) \rtimes G$ is a factor. A group action $G \actson (X,\mu)$ is called W$^*$-superrigid if $G \actson (X,\mu)$ can be entirely recovered from this group measure space construction $L^\infty(X) \rtimes G$. This is a strictly stronger property than OE superrigidity and both properties coincide if one can prove that $M$ has a unique (group measure space) Cartan subalgebra, see e.g.\ \cite{PV09}. When dealing with actions that are not measure preserving, and especially with actions of type III$_0$ the same nuances as with OE superrigidity appear and we get the natural definitions of W$^*$-superrigidity (v1) and (v2).

For none of the concrete actions in Corollary \ref{cor.main-example} and Theorem \ref{thm.main-no-go-theorem}, it is known whether the crossed product has a unique (group measure space) Cartan subalgebra, up to unitary conjugacy. Nevertheless, repeating the construction of \cite[Proposition D]{Vae13}, we obtain ad hoc examples of group actions that are W$^*$-superrigid (v2), of type III$_0$, with a prescribed associated flow. We explain this in Remark \ref{rem.W-star-superrigid}.

\section{Preliminaries}\label{sec.prelim}

Recall that an action of a countable group $G$ on a standard probability space $(X,\mu)$ is called nonsingular if it preserves Borel sets of measure zero. We write $(g \cdot \mu)(\cU) = \mu(g^{-1} \cdot \cU)$ and consider the Radon-Nikodym derivatives $d(g \cdot \mu) / d\mu$, which are well defined a.e. Given a nonsingular action $G \actson (X,\mu)$ of a countable group $G$ on a standard probability space $(X,\mu)$, we consider the associated \emph{Maharam extension}
\begin{equation}\label{eq.maharam}
G \actson X \times \R : g \cdot (x,s) = (g \cdot x, \om(g,x) + s) \; ,
\end{equation}
where $\om(g,x) = \log (d(g^{-1} \cdot \mu)/d\mu)(x)$ is the logarithm of the Radon-Nikodym $1$-cocycle. We may and always will assume that $\om$ is a strict $1$-cocycle, meaning that the cocycle identity holds everywhere. We equip $X \times \R$ with the $G$-invariant $\sigma$-finite measure $d\mu(x) \times e^{-s}\,ds$. One considers the measure scaling action
\begin{equation}\label{eq.commuting-scaling}
\R \actson X \times \R : t \cdot (x,s) = (x,t+s) \; ,
\end{equation}
which commutes with the Maharam extension $G \actson X \times \R$. Denote by $(Y,\eta)$ the space of ergodic components of $G \actson X \times \R$, together with the nonsingular factor map $\pi : X \times \R \to Y$.
Since the actions of $G$ and $\R$ on $X \times \R$ commute, there is an essentially unique nonsingular action $\R \actson (Y,\eta)$ such that for all $t \in \R$, we have that $\pi(t \cdot (x,s)) = t \cdot \pi(x,s)$ for a.e.\ $(x,s) \in X \times \R$. The action $\R \actson (Y,\eta)$ is Krieger's \emph{associated flow} of the action $G \actson (X,\mu)$.

By \cite[Proposition B.5]{Zim84}, after discarding from $(X,\mu)$ a $G$-invariant Borel null set, we may assume that the factor map $\pi : X \times \R \to Y$ is strictly $G$-invariant and $\R$-equivariant, i.e.\ $\pi(g \cdot (x,s)) = \pi(x,s)$ and $\pi(x,s+t) = t \cdot \pi(x,s)$ for all $g \in G$, $t \in \R$ and $(x,s) \in X \times \R$. Writing $\psi(x) = \pi(x,0)$, we have found a Borel map $\psi : X \to Y$ satisfying
\begin{equation}\label{eq.choice-psi}
\pi(x,s) = s \cdot \psi(x) \quad\text{and}\quad \psi(g \cdot x) = (-\om(g,x)) \cdot \psi(x)
\end{equation}
for all $x \in X$, $s \in \R$, $g \in G$.

Let $G \actson (X,\mu)$ be an essentially free, ergodic, nonsingular action of a countable group $G$ on a nonatomic standard probability space $(X,\mu)$. Recall that the \emph{type} of this action is defined as follows: if there exists a $G$-invariant probability measure $\nu \sim \mu$, the action is of type II$_1$~; if there exists a $G$-invariant infinite measure $\nu \sim \mu$, the action is of type II$_\infty$~; in all other cases, the action is of type III. Also recall that $G \actson (X,\mu)$ is of type II$_1$ or II$_\infty$ if and only if the associated flow is isomorphic with the translation action $\R \actson \R$. When the associated flow is not the translation action, there are three possibilities: if $Y$ is reduced to one point, the action is said to be of type III$_1$~; if the associated flow is isomorphic with the periodic flow $\R \actson \R / \Z \log \lambda$ with $0 < \lambda < 1$, the action is said to be of type III$_\lambda$~; finally, when the associated flow is properly ergodic, the action is said to be of type III$_0$.

Two nonsingular actions $G \actson (X,\mu)$ and $\Lambda \actson (Z,\zeta)$ are said to be \emph{conjugate} if there exists an isomorphism of groups $\delta : G \to \Lambda$ and a nonsingular isomorphism $\Delta : (X,\mu) \to (Z,\zeta)$ such that $\Delta(g \cdot x) = \delta(g) \cdot \Delta(x)$ for all $g \in G$ and a.e.\ $x \in X$. Two nonsingular actions $G \actson (X,\mu)$ and $G \actson (Z,\zeta)$ of the same group $G$ are said to be \emph{isomorphic} if there exists a nonsingular isomorphism $\Delta : (X,\mu) \to (Z,\zeta)$ such that $\Delta(g \cdot x) = g \cdot \Delta(x)$ for all $g \in G$ and a.e.\ $x \in X$.

Two essentially free, ergodic, nonsingular actions $G \actson (X,\mu)$ and $\Lambda \actson (Z,\zeta)$ are called \emph{stably orbit equivalent} if there exist nonnegligible Borel sets $\cU \subset X$, $\cV \subset Z$ and a nonsingular isomorphism $\Delta : \cU \to \cV$ such that $\Delta(\cU \cap G \cdot x) = \cV \cap \Lambda \cdot \Delta(x)$ for a.e.\ $x \in \cU$. The actions are called \emph{orbit equivalent} if we may choose $\cU = X$ and $\cV = Z$. When the actions are both of type II$_\infty$ or type III, stable orbit equivalence is the same as orbit equivalence. Recall that the associated flow is invariant under stable orbit equivalence.

We say that a nonsingular action $G \actson (X,\mu)$ is \emph{induced} if there exists a proper subgroup $G_0 < G$ and a $G_0$-invariant Borel set $X_0 \subset X$ such that the sets $(g \cdot X_0)_{g \in G/G_0}$ are disjoint and $\mu(X \setminus G \cdot X_0) = 0$. We then say that $G \actson X$ is induced from $G_0 \actson X_0$. Given any nonsingular action $G_0 \actson (X_0,\mu_0)$ and a larger countable group $G$ containing $G_0$, there is, up to isomorphism, a unique nonsingular action $G \actson X$ that is induced from $G_0 \actson X_0$. Note that by construction, if $G \actson X$ is induced from $G_0 \actson X_0$, then $G \actson X$ and $G_0 \actson X_0$ are stably orbit equivalent. For later reference, we record the following lemma.

\begin{lemma}\label{lem.not-induced}
An ergodic nonsingular action $G \actson (X,\mu)$ is not induced if and only if for every action $G \actson I$ of $G$ on a countable set $I$, every $G$-equivariant Borel map $X \to I$ is constant a.e.
\end{lemma}
\begin{proof}
If $G \actson (X,\mu)$ is not induced, $G \actson I$ and $\vphi : X \to I$ is $G$-equivariant, we can take $i_0 \in I$ such that $X_0 = \{x \in X \mid \vphi(x) = i_0\}$ is nonnegligible. Defining $G_0 = \{g \in G \mid g \cdot i_0 = i_0\}$, it follows that $X_0$ is $G_0$-invariant and that $(g \cdot X_0)_{g \in G/G_0}$ is a partition of $X$, up to measure zero. Since $G \actson X$ is not induced, it follows that $G_0 = G$ and that $X_0 = X$, up to measure zero. This means that $\vphi$ is essentially constant.

If $G \actson X$ is induced from $G_0 \actson X_0$, the map $x \mapsto gG_0$ for $x \in g \cdot X_0$ is a $G$-equivariant map $X \to G/G_0$ that is not essentially constant.
\end{proof}

An essentially free, nonsingular action $\Sigma \actson (X,\mu)$ is said to \emph{admit a fundamental domain} if there exists a Borel set $\cU \subset X$ such that all $g \cdot \cU$, $g \in \Sigma$, are disjoint and $\mu(X \setminus \Sigma \cdot \cU) = 0$. In that case, the quotient $X/\Sigma$ is a well defined standard measure space and identified with $(\cU,\mu)$.

\section{\boldmath Adjoint flows and type III$_0$ actions with prescribed associated flow}

We say that a flow $\R \actson (Z,\zeta)$ \emph{scales} the $\si$-finite measure $\zeta$ if $t \cdot \zeta = e^t \, \zeta$ for all $t \in \R$.

\begin{proposition}\label{prop.unique-char-adj}
Let $\R \actson (Y,\eta)$ be an ergodic flow. Up to isomorphism, there is a unique nonsingular ergodic action $\R^2 \actson (Z,\zeta)$ of $\R^2$ on a standard, $\si$-finite measure space $(Z,\zeta)$ such that the actions of both $\R \times \{0\}$ and $\{0\} \times \R$ scale the measure $\zeta$ and such that $\R \actson Z/(\{0\} \times \R)$ is isomorphic with $\R \actson Y$.
\end{proposition}
\begin{proof}
Denote by $\om : \R \times Y \to \R$ the logarithm of the Radon-Nikodym cocycle. Define the measure $\gamma$ on $\R$ by $d\gamma(t) = e^{-t} \, dt$. Define $(Z,\zeta) = (Y \times \R,\eta \times \gamma)$ and define the action
$$\R^2 \actson (Z,\zeta) : (t,r) \cdot (y,s) = (t \cdot y, \om(t,y) + t + r + s) \; .$$
Both the actions by $(t,0)$ and by $(0,r)$ scale the measure $\zeta$. By construction, $Z/(\{0\} \times \R) = Y$.

Now assume that $\R^2 \actson (Z',\zeta')$ is a nonsingular ergodic action such that the actions of both $\R \times \{0\}$ and $\{0\} \times \R$ scale the measure $\zeta'$ and such that $\R \actson Z'/(\{0\} \times \R)$ is isomorphic with $\R \actson Y$. We prove that $\R^2 \actson (Z',\zeta')$ is isomorphic with $\R^2 \actson (Z,\zeta)$.

Since the action of $\{0\} \times \R$ scales the measure $\zeta'$ and $\R \actson Z'/(\{0\} \times \R)$ is isomorphic with $\R \actson Y$, we find a $\si$-finite measure $\mu' \sim \mu$ on $Y$ and a measure preserving isomorphism $\Delta : (Y \times \R,\mu' \times \gamma) \to (Z',\zeta')$ such that for all $(t,r) \in \R^2$,
$$\Delta(t \cdot y, \zeta(t,y) + r + s) = (t,r) \cdot \Delta(y,s) \quad\text{for a.e.\ $(y,s) \in Y \times \R$,}$$
where $\zeta : \R \times Y \to \R$ is a $1$-cocycle. Precomposing $\Delta$ with the measure preserving map
$$(Y \times \R,\mu \times \gamma) \to (Y \times \R,\mu' \times \gamma) : (y,s) \mapsto (y, \log(d\mu' / d\mu)(y) + s)$$
and replacing $\zeta$ by a cohomologous $1$-cocycle, we may assume that $\mu' = \mu$. Expressing that the action of $\R \times \{0\}$ scales the measure $\mu \times \gamma$ gives us that $\zeta(t,y) = \om(t,y) + t$. So we have found the required isomorphism.
\end{proof}

Given the uniqueness of $\R^2 \actson (Z,\zeta)$ in Proposition \ref{prop.unique-char-adj}, we get the following well defined notion of an adjoint flow and we automatically have that this adjoint is an involutive operation: the adjoint of $\alh$ is isomorphic with $\al$.

\begin{definition}\label{def.adjoint-flow}
Given an ergodic flow $\R \actson^\al (Y,\eta)$, the adjoint flow $\R \actson^{\alh} (\Yh,\etah)$ is defined as the ergodic flow $\R \actson Z/(\R \times \{0\})$, where $\R^2 \actson Z$ is the unique action given by Proposition \ref{prop.unique-char-adj}.
\end{definition}

Note that we can also define the adjoint flow $\alh$ more concretely. Denoting by $\om : \R \times Y \to \R$ the logarithm of the Radon-Nikodym cocycle of an ergodic flow $\R \actson^\al Y$, we consider the quotient $(Y \times \R)/\R$, where $\R$ is acting by $t \cdot (y,s) = (t \cdot y, t + \om(t,y) + s)$. On this quotient, we let $\R$ act by translation in the second variable. This is the adjoint flow $\alh$.

From this concrete description, it immediately follows that $\al \cong \alh$ whenever the flow $\al$ admits a finite or $\si$-finite equivalent $\R$-invariant measure. In general, $\al$ need not be isomorphic with $\alh$, as the following example shows.

\begin{example}\label{def.example-adjoint-flow}
Let $\R \actson^\al (Y,\eta)$ be the ergodic flow given as the induction of an ergodic, type III$_1$ action $\Z \actson^{\al_0} (Y_0,\eta_0)$. We prove that the adjoint flow $\alh$ is not isomorphic with $\al$.

Recall that the induced flow $\al$ is defined as follows. Consider the action $\R \times \Z \actson \R \times Y_0 : (t,n) \cdot (s,y) = (t-n+s,n \cdot y)$. Then, $\al$ is defined as the action of $\R$ on $(\R \times Y_0)/(\{0\} \times \Z)$. To determine the adjoint flow $\alh$, denote by $\om : \Z \times Y_0 \to \R$ the logarithm of the Radon-Nikodym cocycle for $\al_0$. Denote by $\lambda$ the Lebesgue measure on $\R$ and define the measure $\gamma$ such that $(d\gamma/d\lambda)(t) = e^{-t}$. Then consider the action
$$\R \times \Z \times \R \actson (\R \times Y_0 \times \R,\lambda \times \eta_0 \times \gamma) : (t,n,r) \cdot (s,y,s') = (t - n + s, n \cdot y, t+\om(n,y) + r + s') \; .$$
The action of $\{0\} \times \Z \times \{0\}$ is measure preserving. The actions of $\R \times \{(0,0)\}$ and $\{(0,0)\} \times \R$ are measure scaling. By construction, $\al$ is given by $\R \actson (\R \times Y_0 \times \R)/(\{0\} \times \Z \times \R)$. We conclude that the adjoint flow $\alh$ is given by $\R \actson (\R \times Y_0 \times \R)/(\R \times \Z \times \{0\})$.

By construction, the flow $\al$ comes with an $\R$-equivariant map $Y \to \R/\Z$. We prove that such an $\R$-equivariant map does not exist for the adjoint flow $\alh$. Assuming the contrary, we find a map $\theta : \R \times Y_0 \times \R \to \R/\Z$ that is invariant for the action of $\R \times \Z \times \{0\}$ and that is equivariant for the action of $\{(0,0)\} \times \R$. By the invariance under $\R \times \{(0,0)\}$ and the equivariance under $\{(0,0)\} \times \R$, the map $\theta$ must be of the form $\theta(s,y,s') = -s + s' + \vphi(y)$ where $\vphi : Y_0 \to \R/\Z$. The invariance under $\{0\} \times \Z \times \{0\}$ then says that
$$\om(n,y) + \vphi(n \cdot y) = \vphi(y) \quad\text{for all $n \in \Z$ and a.e.\ $y \in Y_0$.}$$
This means that the map $Y_0 \times \R \to \R/\Z : (y,s) \mapsto s + \vphi(y)$ is invariant under the action $\Z \actson Y_0 \times \R : n \cdot (y,s) = (n\cdot y, \om(n,y) + s)$. This action is ergodic because $\Z \actson Y_0$ is assumed to be of type III$_1$. So, the map $(y,s) \mapsto s + \vphi(y)$ is essentially constant, which is absurd. So we have proven that $\alh$ is not isomorphic with $\al$.
\end{example}

\begin{proposition}\label{prop.prescribed-flow}
Let $G \actson (X,\mu)$ be a nonsingular ergodic action of type III$_1$, with logarithm of the Radon-Nikodym cocycle $\om : G \times X \to \R$. Let $\R \actson^\al (Y,\eta)$ be any ergodic flow. Then the action
\begin{equation}\label{eq.my-good-action}
G \actson (X \times Y, \mu \times \eta) : g \cdot (x,y) = (g \cdot x, \al_{\om(g,x)}(y))
\end{equation}
is ergodic and has associated flow $\alh$.
\end{proposition}
\begin{proof}
Let $\R^2 \actson (Z,\zeta)$ be the unique action given by Proposition \ref{prop.unique-char-adj}, associated with the ergodic flow $\R \actson (Y,\eta)$. Consider the action
\begin{equation}\label{eq.my-Maharam}
G \times \R \actson (X \times Z,\mu \times \zeta) : (g,r) \cdot (x,z) = (g\cdot x, (\om(g,x),r) \cdot z) \; .
\end{equation}
The action of $\{e\} \times \R$ scales the measure, while the action of $G \times \{0\}$ is measure preserving. By definition, the action of $G$ on $(X \times Z)/(\{e\} \times \R)$ is isomorphic with $G \actson X \times Y$. It thus follows that the action in \eqref{eq.my-Maharam} is the Maharam extension of $G \actson X \times Y$ together with its measure scaling action of $\R$.

By the uniqueness of $\R^2 \actson (Z,\zeta)$, we may as well identify $(Z,\zeta) = (\R \times \Yh, \gamma \times \etah)$ with
$$(t,r) \cdot (s,\yh) = (t + r+\beh(r,\yh) + s , r \cdot \yh) \; ,$$
where $\beh : \R \times \Yh \to \R$ is the logarithm of the Radon-Nikodym cocycle for the adjoint flow $\R \actson^{\alh} \Yh$. Then the action in \eqref{eq.my-Maharam} becomes
$$G \times \R \actson (X \times \R \times \Yh, \mu \times \gamma \times \etah) : (g,r) \cdot (x,s,\yh) = (g \cdot x, \om(g,x) + \beh(r,\yh) + r + s, r \cdot \yh) \; .$$
Since $G \actson (X,\mu)$ is ergodic and of type III$_1$, the Maharam extension $G \actson X \times \R$ is ergodic. It follows that the $G$-invariant functions on $X \times \R \times \Yh$ are the functions that only depend on the $\Yh$-variable. Since $\R \actson \Yh$ is ergodic, we conclude that the action in \eqref{eq.my-Maharam} is ergodic. We have proven that the action in \eqref{eq.my-good-action} is ergodic and that its associated flow is identified with $\alh : \R \actson \Yh$.
\end{proof}

\section{Versions of OE superrigidity in the type III setting}\label{sec.versions-OE-superrigidity}

\begin{definition}\label{def.simple-action}
We say that a free, ergodic, nonsingular action $G \actson (X,\mu)$ of a countable group $G$ is \emph{simple} if the action is not induced and if there are no nontrivial normal subgroups $\Sigma \lhd G$ for which $\Sigma \actson (X,\mu)$ admits a fundamental domain.
\end{definition}

The motivation for this ad hoc notion of simplicity is the following. When $G \actson (X,\mu)$ is induced from $G_0 \actson X_0$, we have a canonical stable orbit equivalence between $G \actson X$ and $G_0 \actson X_0$. When $\Sigma \lhd G$ is a normal subgroup such that $\Sigma \actson (X,\mu)$ admits a fundamental domain, we have a canonical stable orbit equivalence between $G \actson X$ and $G/\Sigma \actson X/\Sigma$. So when $G \actson (X,\mu)$ is not simple, there always is a certain absence of OE-superrigidity and describing all stably orbit equivalent actions is necessarily cumbersome. For this reason, we mainly restrict ourselves to simple actions in this paper.

Recall that a nonsingular action $G \actson (X,\mu)$ of a countable group $G$ on a standard probability space $(X,\mu)$ is called \emph{cocycle superrigid with countable target groups} if every $1$-cocycle $\Om : G \times X \to \Lambda$ with values in a countable group $\Lambda$ is cohomologous to a group homomorphism $\delta : G \to \Lambda$, viewed as a $1$-cocycle that is independent of the $X$-variable.

Recall from the introduction the two versions (v1) and (v2) of OE superrigidity. For simple actions, version (v1) of OE superrigidity turns out to be equivalent with cocycle superrigidity with countable targets.

\begin{proposition}\label{prop.OE-v1-versus-cocycle}
Let $G \actson (X,\mu)$ be any free, ergodic, nonsingular, simple action. Then $G \actson (X,\mu)$ satisfies OE superrigidity (v1) if and only if $G \actson (X,\mu)$ is cocycle superrigid with countable targets.
\end{proposition}

\begin{proof}
The implication from cocycle superrigidity to OE superrigidity was first proven, in a pmp setting, in \cite[Proposition 4.2.11]{Zim84}. The version that we need is literally proven in \cite[Lemma 2.4]{DV21}.

Conversely, assume that $\Om : G \times X \to \Lambda$ is a $1$-cocycle with values in a countable group $\Lambda$. Consider the free, nonsingular, ergodic action
$$G \times \Lambda \actson X \times \Lambda : (g,a) \cdot (x,b) = (g \cdot x, \Om(g,x) b a^{-1}) \; .$$
By construction, this action is stably orbit equivalent with $G \actson X$. Assume that this action is conjugate to an induction of $G \actson X$. We have to prove that $\Om$ is cohomologous to a group homomorphism.

Take an injective group homomorphism $\delta : G \to G \times \Lambda: \delta(g) = (\delta_1(g),\delta_2(g))$ and a measure space isomorphism $\Delta : X \to Z \subset X \times \Lambda$ such that $G \times \Lambda \actson X \times \Lambda$ is induced from $\delta(G) \actson Z$ and $\Delta$ is a conjugacy w.r.t.\ $\delta$.

Since the action of $\Lambda$ on $X \times \Lambda$ admits a fundamental domain, the same is true for the action $\Ker \delta_1 \actson X$. Since $G \actson X$ is simple, we find that $\delta_1$ is faithful. Since $\delta(G) \subset \delta_1(G) \times \Lambda$ and since $G \times \Lambda \actson X \times \Lambda$ is induced from $\delta(G) \actson Z$, we find a fortiori that $G \times \Lambda \actson X \times \Lambda$ is induced from $\delta_1(G) \times \Lambda \actson Z_1$ with $Z_0 \subset Z_1$. Since $Z_1$ is $\Lambda$-invariant, we find that $Z_1 = X_0 \times \Lambda$ and conclude that $G \actson X$ is induced from $\delta_1(G) \actson X_0$. Since $G \actson X$ is simple, we conclude that $\delta_1(G) = G$. So, $\delta_1$ is an automorphism of the group $G$.

Define the group homomorphism $\gamma : G \to \Lambda : \gamma = \delta_2 \circ \delta_1^{-1}$. We have $\delta(g) = (\delta_1(g),\gamma(\delta_1(g)))$ and the map
$$\psi : (G \times \Lambda) / \delta(G) \to \Lambda : (g,k)\delta(G) \mapsto \gamma(g) k^{-1}$$
is a bijection satisfying $\psi((g,a) \cdot i) = \gamma(g) \psi(i) a^{-1}$ for all $(g,a) \in G \times \Lambda$ and $i \in (G \times \Lambda)/\delta(G)$.

Since $G \times \Lambda \actson X \times \Lambda$ is induced from $\delta(G) \actson Z$, we find a $(G \times \Lambda)$-equivariant map from $X \times \Lambda$ to $(G \times \Lambda) / \delta(G)$. We denote by $\theta$ its composition with $\psi$. By $\Lambda$-equivariance, we get that $\theta(x,a) = \theta_0(x) a$, where $\theta_0 : X \to \Lambda$ is a Borel map. Expressing the $G$-equivariance gives us that
$$\Om(g,x) = \theta_0(g \cdot x)^{-1} \, \gamma(g) \, \theta_0(x) \; ,$$
so that $\Om$ is cohomologous to a group homomorphism.
\end{proof}


As an essential ingredient to prove the first part of Theorem \ref{thm.main-no-go-theorem}, as well as to prove Corollary \ref{cor.main-example}, we need to establish cocycle superrigidity for linear and for affine actions on $\R^n$. We first recall the notion of essential cocycle superrigidity introduced in \cite[Definition B]{DV21}. We only formulate the version for connected Lie groups, which is the one that we need in this paper.

\begin{definition}[{\cite[Definition B]{DV21}}]
A countable dense subgroup $\Gamma < G$ of a connected Lie group $G$, with universal cover $\pi : \Gtil \to G$ and $\Gammatil = \pi^{-1}(\Gamma)$, is said to be essentially cocycle superrigid with countable targets if for every $1$-cocycle $\om : \Gamma \times G \to \Lambda$ of the translation action $\Gamma \actson G$ with values in a countable group $\Lambda$, the lifted $1$-cocycle $\omtil : \Gammatil \times \Gtil \to \Lambda : \omtil = \om \circ (\pi \times \pi)$ is cohomologous to a group homomorphism $\delta : \Gammatil \to \Lambda$.
\end{definition}

In \cite[Propositions 4.1 and 4.2]{DV21}, it was proven that for every integer $n \geq 3$, nonempty set of prime numbers $\cS$ and real algebraic number field $\Q \varsubsetneq K \subset \R$ with ring of integers $\cO_K$, the dense subgroups $\SL(n,\Z[\cS^{-1}])$ and $\SL(n,\cO_K)$ of $\SL(n,\R)$ are essentially cocycle superrigid with countable targets. We prove the same result for a much larger family of dense subgroups of $\SL(n,\R)$ and also for their corresponding subgroups of $\SL(n,\R) \ltimes \R^n$.

For every integer $n \geq 2$ and commutative ring $\cA$, we denote by $\SL(n,\cA)$ and $\GL(n,\cA)$ the groups of $n \times n$ matrices with entries in $\cA$ and determinant resp.\ $1$, belonging to $\cA^*$. Whenever $1 \leq i,j \leq n$ with $i \neq j$ and $a \in \cA$, we denote by $e_{ij}(a) \in \SL(n,\cA)$ the elementary matrix with $1$'s on the diagonal, $a$ in position $ij$ and $0$'s elsewhere. We denote by $E(n,\cA) \subset \SL(n,\cA)$ the subgroup generated by the elementary matrices. By Suslin's theorem (see \cite[Theorem 1.2.13]{HOM89}), for all $n \geq 3$, $E(n,\cA)$ is a normal subgroup of $\GL(n,\cA)$. For several rings $\cA$, it is known that $E(n,\cA) = \SL(n,\cA)$ (see Example \ref{ex.allerlei-voorbeelden} for references).

\begin{theorem}\label{thm.essential-cocycle-superrigidity}
Let $\cA \subset \R$ be any countable subring containing an algebraic number that does not belong to $\Z$. Let $n \geq 3$ and $E(n,\cA) < \Gamma < \SL(n,\cA)$ any intermediate subgroup.
\begin{enumlist}
\item The dense subgroup $\Gamma < \SL(n,\R)$ is essentially cocycle superrigid with countable targets.
\item The dense subgroup $\Gamma \ltimes \cA^n < \SL(n,\R) \ltimes \R^n$ is essentially cocycle superrigid with countable targets.
\end{enumlist}
\end{theorem}

Before proving Theorem \ref{thm.essential-cocycle-superrigidity}, we need two elementary lemmas, which are essentially contained in \cite[Lemma 5.1]{Ioa14}.

\begin{lemma}\label{lem.dense-not-induced}
Let $G$ be a connected locally compact second countable group and let $\Gamma < G$ be a countable dense subgroup. Then the translation action $\Gamma \actson G$ is ergodic and not induced.
\end{lemma}
\begin{proof}
Let $\Gamma_0 < \Gamma$ be a subgroup and $\pi : G \to \Gamma/\Gamma_0$ a Borel map such that $\pi(g h) = g \pi(h)$ for all $g \in \Gamma$ and a.e.\ $h \in G$. We have to prove that $\pi$ is essentially constant. Since $\Gamma/\Gamma_0$ is countable, the set $\cU = \{(h,k) \in G \times G \mid \pi(h) = \pi(k)\}$ is nonnegligible.

Since $\pi(g h) = g \pi(h)$ when $g \in \Gamma$, the set $\cU$ is essentially invariant under the diagonal translation action $\Gamma \actson G \times G$. By continuity, $\cU$ is essentially invariant under the diagonal translation action of $G$. We thus find a nonnegligible Borel set $\cV \subset G$ such that, up to measure zero, $\cU = \{(h,h k) \mid h \in G, k \in \cV\}$.

Define $G_0 = \{k \in G \mid \pi(h) = \pi(hk) \;\text{for a.e.\ $h \in G$}\}$. Since $\cV$ is nonnegligible, also $G_0$ is nonnegligible. By definition, $G_0$ is a subgroup of $G$. So, $G_0$ must be an open subgroup of $G$. Since $G$ is connected, we conclude that $G_0 = G$. This means that $\pi$ is essentially constant.
\end{proof}

Also the following lemma is essentially contained in \cite[Lemma 5.1]{Ioa14} and allows to extend cocycle superrigidity from a subgroup $\Gamma_0 < \Gamma$ to its normalizer $\cN_\Gamma(\Gamma_0)$. Such a result goes back to \cite[Proposition 3.6]{Pop05}.


\begin{lemma}\label{lem.further-trivialize}
Let $\Gamma \actson (X,\mu)$ be a free, ergodic, nonsingular action and let $\Om : \Gamma \times X \to \Lambda$ be a $1$-cocycle with values in a countable group $\Lambda$. Let $\Gamma_0 < \Gamma$ be a subgroup and let $(X_y,\mu_y)_{y \in (Y,\eta)}$ be the ergodic decomposition of the action $\Gamma_0 \actson (X,\mu)$, with corresponding factor map $\pi : (X,\mu) \to (Y,\eta)$. Assume that for $\eta$-a.e.\ $y \in Y$, the action $\Gamma_0 \actson (X_y,\mu_y)$ is not induced.

If for every $h \in \Gamma_0$, the function $x \mapsto \Om(h,x)$ factors through $\pi$, then $x \mapsto \Om(g,x)$ factors through $\pi$ for every $g \in \cN_\Gamma(\Gamma_0)$.
\end{lemma}
\begin{proof}
Fix $g \in \cN_\Gamma(\Gamma_0)$ and denote by $\al : \Gamma_0 \to \Gamma_0$ the automorphism $\al(h) = g h g^{-1}$. Since for every $h \in \Gamma_0$, the function $x \mapsto \Om(h,x)$ factors through $\pi$, we find a measurable family $(\delta_y)_{y \in Y}$ of group homomorphisms $\delta_y : \Gamma_0 \to \Lambda$ such that $\Om(h,x) = \delta_{\pi(x)}(h)$ for all $h \in \Gamma_0$ and a.e.\ $x \in X$. Applying the $1$-cocycle relation to $g h = \al(h) g$ gives us that $\Om(g,h \cdot x) \, \delta_{\pi(x)}(h) = \delta_{\pi(g \cdot x)}(\al(h)) \, \Om(g,x)$ for a.e.\ $x \in X$.

For every $g \in \cN_\Gamma(\Gamma_0)$, the map $x \mapsto \pi(g \cdot x)$ is $\Gamma_0$-invariant. We can thus define the nonsingular action $\cN_\Gamma(\Gamma_0) \actson (Y,\eta)$ such that $\pi(g \cdot x) = g \cdot \pi(x)$ for all $g \in \cN_\Gamma(\Gamma_0)$ and a.e.\ $x \in X$. We conclude that for a.e.\ $y \in Y$, we have
\begin{equation}\label{eq.equivariance}
\Om(g,h \cdot x) = \delta_{g \cdot y}(\al(h)) \, \Om(g,x) \, \delta_y(h)^{-1} \quad\text{for $\mu_y$-a.e.\ $x \in X_y$.}
\end{equation}
The action $\Gamma_0 \actson (X_y,\mu_y)$ is ergodic and not induced. Also \eqref{eq.equivariance} is saying that the map $X_y \to \Lambda : x \mapsto \Om(g,x)$ is $\Gamma_0$-equivariant, where $\Gamma_0$ is acting on $\Lambda$ by $h \cdot \lambda = \delta_{g \cdot y}(\al(h)) \, \lambda \, \delta_y(h)^{-1}$. By Lemma \ref{lem.not-induced}, the map $x \mapsto \Om(g,x)$ is essentially constant on $(X_y,\mu_y)$. Since this holds for a.e.\ $y \in Y$, we have proven that $x \mapsto \Om(g,x)$ factors through $\pi$.
\end{proof}

Note that Lemma \ref{lem.dense-not-induced} also implies the following result, which we will use in combination with Lemma \ref{lem.further-trivialize}~: if $G$ is a locally compact second countable group and $\Gamma < G$ is a countable subgroup whose closure $H = \overline{\Gamma}$ is connected, then the map $\pi : G \to H \backslash G : \pi(g) = H g$ realizes the ergodic decomposition of the translation action $\Gamma \actson G$ and, for every $H g \in H \backslash G$, the action of $\Gamma$ on $\pi^{-1}(Hg)$ is isomorphic with the translation action $\Gamma \actson H$ and thus not induced by Lemma \ref{lem.dense-not-induced}.

\begin{proof}[{Proof of Theorem \ref{thm.essential-cocycle-superrigidity}}]
By our assumption, the ring $\cA$ contains a rational number $q \in \Q \setminus \Z$ or an irrational algebraic number $\al$. In the first case, by taking a multiple of $q$, we find a prime $p$ and a positive integer $N \in \N \setminus \{0\}$ with $p \nmid N$ such that $N p^{-1} \in \cA$. Since $p \nmid N$, we can take integers $a,b \in \Z$ such that $aN + bp = 1$, so that $p^{-1} = a Np^{-1} + b$. We conclude that $p^{-1} \in \cA$ and denote $\cA_0 = \Z[p^{-1}]$. In the second case, we define the algebraic number field $K = \Q(\al)$. Let $d \geq 2$ be the degree of the minimal polynomial of $\al$. The ring $\cO_K$ of integers of $K$ is a finitely generated $\Z$-module that is contained in $K$, which has $\{1,\al,\ldots,\al^{d-1}\}$ as a $\Q$-vector space basis. We can thus take a positive integer $N \in \N \setminus \{0\}$ such that $N \cO_K \subset \Z[\al] \subset \cA$. We denote $\cA_0 = \Z + N \cO_K$.

We start by proving that $E(n,\cA_0) < \SL(n,\R)$ is essentially cocycle superrigid with countable targets. When $\cA_0 = \Z[p^{-1}]$, we know by \cite[Theorem 4.3.9]{HOM89} that $E(n,\cA_0) = \SL(n,\cA_0)$ and we know from \cite[Proposition 4.1]{DV21} that $\SL(n,\cA_0) < \SL(n,\R)$ is essentially cocycle superrigid with countable targets.

Next consider the case where $\cA_0 = \Z + N \cO_K$. The ring $\cA_1 = \cO_K / (N \cO_K)$ is finite and the kernel of the canonical homomorphism $\SL(n,\cO_K) \to \SL(n,\cA_1)$ is contained in $\SL(n,\cA_0)$. Thus, $\SL(n,\cA_0) < \SL(n,\cO_K)$ has finite index. Using the real and complex embeddings of $K$, the group $\SL(n,\cO_K)$ is an irreducible lattice in a product of copies of $\SL(n,\R)$ and $\SL(n,\C)$, see e.g.\ \cite[Theorems 5.7 and 7.12]{PR94}. The groups $\SL(n,\R)$ and $\SL(n,\C)$ both have property~(T), see e.g.\ \cite[Theorem 1.4.15]{BHV08}. Then also the finite index subgroup $\SL(n,\cA_0)$ of $\SL(n,\cO_K)$ is such an irreducible lattice. In particular, $\SL(n,\cA_0)$ satisfies Margulis' normal subgroup theorem. Since $E(n,\cA_0)$ is an infinite normal subgroup of $\SL(n,\cA_0)$, we conclude that $E(n,\cA_0) < \SL(n,\cA_0)$ has finite index. So also $E(n,\cA_0)$ is an irreducible lattice in a product of copies of $\SL(n,\R)$ and $\SL(n,\C)$. By part 2 of \cite[Theorem C]{DV21} (and this is essentially \cite[Theorem B]{Ioa14}), we get that $E(n,\cA_0) < \SL(n,\R)$ is essentially cocycle superrigid with countable targets.

We also need the following observation: given a countable dense subgroup $\Gamma$ of a connected Lie group $G$ with universal cover $\pi : \Gtil \to G$, the subgroup $\Gamma < G$ is essentially cocycle superrigid with countable targets if and only if the translation action of the dense subgroup $\Gammatil = \pi^{-1}(\Gamma)$ on $\Gtil$ is cocycle superrigid with countable targets. One implication is obvious. So assume that $\Gamma < G$ is essentially cocycle superrigid with countable targets and that $\om : \Gammatil \times \Gtil \to \Lambda$ is a $1$-cocycle with values in a countable group $\Lambda$. Since the action of $\Ker \pi$ on $\Gtil$ admits a fundamental domain, $\om$ is cohomologous with a $1$-cocycle $\om_1$ satisfying $\om_1(g,x) = e$ for all $g \in \Ker \pi$ and a.e.\ $x \in \Gtil$. This means that $\om_1(g,x) = \om_2(\pi(g),\pi(x))$ for all $g \in \Gammatil$ and a.e.\ $x \in \Gtil$, where $\om_2 : \Gamma \times G \to \Lambda$ is a $1$-cocycle. By our assumption, $\om_1$, and hence also $\om$, is cohomologous with a group homomorphism from $\Gammatil$ to $\Lambda$.

We then prove statement 1. Take $E(n,\cA) < \Gamma < \SL(n,\cA)$. We write $G = \SL(n,\R)$ and denote by $\pi : \Gtil \to G$ its universal cover. Since $n \geq 3$, $\Ker \pi$ is a central subgroup of order $2$ in $\Gtil$. Define $\Gammatil = \pi^{-1}(\Gamma)$. Let $\om : \Gamma \times G \to \Lambda$ be a $1$-cocycle, with lift $\omtil : \Gammatil \times \Gtil \to \Lambda$. We have to prove that $\omtil$ is cohomologous to a group homomorphism.

Since we have proven that $E(n,\cA_0) < G$ is essentially cocycle superrigid, it follows from the observation above that $\omtil$ is cohomologous with a $1$-cocycle $\gamma : \Gammatil \times \Gtil \to \Lambda$ that has the property that $x \mapsto \gamma(g,x)$ is essentially constant for every $g \in \pi^{-1}(E(n,\cA_0))$.

Write $[n] = \{1,\ldots,n\}$. For every $k \in [n]$ and for every subring $\cB \subset \R$, we define the following subgroups of $E(n,\cB)$: the group $C_k(\cB) \cong \cB^{n-1}$ generated by $\{e_{ik}(b) \mid i \in [n] \setminus \{k\}, b \in \cB\}$~; the group $R_k(\cB) \cong \cB^{n-1}$ generated by $\{e_{kj}(b) \mid j \in [n] \setminus \{k\}, b \in \cB\}$~; and the group $H_k(\cB) \cong E(n-1,\cB)$ generated by $\{e_{ij}(b) \mid i,j \in [n] \setminus \{k\}, i \neq j, b \in \cB\}$. Note that $H_k(\cB)$ normalizes both $C_k(\cB)$ and $R_k(\cB)$. If $\cB \subset \R$ is dense, also $H_k(\cB)$, $C_k(\cB)$ and $R_k(\cB)$ are dense in resp.\ $H_k(\R)$, $C_k(\R)$ and $R_k(\R)$, and the latter are closed subgroups of $G$.

Since $\R$ is simply connected, there is for all $i \neq j$ a unique continuous group homomorphism $\etil_{ij} : \R \to \Gtil$ such that $\pi(\etil_{ij}(t)) = e_{ij}(t)$ for all $t \in \R$. When $i,j \in [n] \setminus \{k\}$ and $s,t \in \R$, the image $\pi([\etil_{ik}(t),\etil_{jk}(s)])$ of the commutator equals the identity element. Thus, $[\etil_{ik}(t),\etil_{jk}(s)] \in \Ker \pi$ for all $s,t \in \R$. By connectedness of $\R^2$, we find that $\etil_{ik}(t)$ commutes with $\etil_{jk}(s)$. There thus is a unique continuous group homomorphism $C_k(\R) \to \Gtil : e_{ik}(t) \mapsto \etil_{ik}(t)$.

For every $k \in [n]$ and subring $\cB \subset \R$, we denote by $\Ctil_k(\cB)$ the subgroup of $\Gtil$ generated by $\{\etil_{ik}(b) \mid i \in [n] \setminus \{k\}, b \in \cB\}$. Note that $\Ctil_k(\R)$ is a connected closed subgroup of $\Gtil$ and that $\pi : \Ctil_k(\R) \to C_k(\R)$ is an isomorphism. We also have that $\Ctil_k(\R)$ is the connected component of the identity in $\pi^{-1}(C_k(\R))$ and that $\Ctil_k(\cB) = \Ctil_k(\R) \cap \pi^{-1}(C_k(\cB))$.

Define $c_k : \Gtil \to \Ctil_k(\R) \backslash \Gtil : c_k(x) = \Ctil_k(\R) x$. Since $\Ctil_k(\cA)$ commutes with $\Ctil_k(\cA_0)$ and since $\Ctil_k(\cA_0)$ is dense in $\Ctil_k(\R)$, it follows from Lemmas \ref{lem.dense-not-induced} and \ref{lem.further-trivialize} that for every $g \in \Ctil_k(\cA)$, the map $x \mapsto \gamma(g,x)$ factors through $c_k$.

Take $g \in \pi^{-1}(H_k(\cA) C_k(\cA))$. Since $g$ normalizes $\pi^{-1}(C_k(\R))$, the element $g$ also normalizes its connected component of the identity $\Ctil_k(\R)$. Since $g$ normalizes as well $\pi^{-1}(C_k(\cA))$, it follows that $g$ normalizes $\Ctil_k(\cA)$. Another application of Lemmas \ref{lem.dense-not-induced} and \ref{lem.further-trivialize} then says that the map $x \mapsto \gamma(g,x)$ factors through $c_k$.

Fix $i \neq j$ and fix $g \in \pi^{-1}(e_{ij}(\cA))$. In the following paragraphs, we prove that $x \mapsto \gamma(g,x)$ is essentially constant.

We have proven that the map $x \mapsto \gamma(g,x)$ factors through $c_k$ for all $k \neq i$. This means that for all $b \neq i$, $a \neq b$ and $t \in \R$, we have $\gamma(g,\etil_{ab}(t)x) = \gamma(g,x)$ for a.e.\ $x \in \Gtil$. Since we can make an analogous reasoning using the subgroups $R_k$, we also get for all $a \neq j$, $b \neq a$ and $t \in \R$ that $\gamma(g,\etil_{ab}(t)x) = \gamma(g,x)$ for a.e.\ $x \in \Gtil$.

Define $\Ttil = \{h \in \Gtil \mid \gamma(g,hx) = \gamma(g,x) \;\;\text{for a.e.\ $x \in \Gtil$}\;\}$. Then $\Ttil$ is a closed subgroup of $\Gtil$. Define $T = \pi(\Ttil)$. By the previous paragraph, $e_{ab}(\R) \subset T$ when $a \neq b$, $a \neq j$ or $b \neq i$. When $a = j$ and $b =i$, we choose $c \in [n] \setminus \{i,j\}$ and note that for all $t \in \R$,
$$e_{ac}(t) \, e_{cb}(1) \, e_{ac}(-t) \, e_{cb}(-1) = e_{ab}(t) \; ,$$
so that again $e_{ab}(\R) \subset T$. It follows that $T = G$. The closed subgroup $\Ttil \subset \Gtil$ thus has at most index $2$, so that $\Ttil \subset \Gtil$ is open. Since $\Gtil$ is connected, it follows that $\Ttil = \Gtil$. This means that $x \mapsto \gamma(g,x)$ is essentially constant.

We have thus proven that $x \mapsto \gamma(g,x)$ is essentially constant for every $g \in \pi^{-1}(E(n,\cA))$. Since $\pi^{-1}(E(n,\cA))$ is a normal subgroup of $\Gammatil$, a final application of Lemmas \ref{lem.dense-not-induced} and \ref{lem.further-trivialize} implies that $x \mapsto \gamma(g,x)$ is essentially constant for every $g \in \Gammatil$. This concludes the proof of statement 1.

To prove statement 2, we still write $G = \SL(n,\R)$ with universal cover $\pi : \Gtil \to G$. Note that $\Gtil \ltimes \R^n$ is the universal cover of $G \ltimes \R^n$. Let $\om : (\Gamma \ltimes \cA^n) \times (G \ltimes \R^n) \to \Lambda$ be a $1$-cocycle, with lift $\omtil : (\Gammatil \ltimes \cA^n) \times (\Gtil \ltimes \R^n) \to \Lambda$.

By the same argument as in the beginning of this proof, $E(n,\cA_0) \ltimes \cA_0^n$ is an irreducible lattice in a product of copies of $\SL(n,\R) \ltimes \R^n$ and $\SL(n,\C) \ltimes \C^n$. By \cite[Corollary 1.4.16]{BHV08}, these groups have property~(T). By part 2 of \cite[Theorem C]{DV21}, it follows that the dense subgroup $E(n,\cA_0) \ltimes \cA_0^n$ of $\SL(n,\R) \ltimes \R^n$ is essentially cocycle superrigid with countable targets.

By the observation at the beginning of the proof, we thus find that $\omtil$ is cohomologous with a $1$-cocycle $\gamma$ such that $x \mapsto \gamma(g,x)$ is essentially constant for every $g \in \pi^{-1}(E(n,\cA_0)) \ltimes \cA_0^n$. Denote by $\theta : \Gtil \ltimes \R^n \to \Gtil$ the natural quotient map. Since $\cA_0^n$ is dense in $\R^n$ and since $\cA^n$ commutes with $\cA_0^n$, it follows from Lemmas \ref{lem.dense-not-induced} and \ref{lem.further-trivialize} that $x \mapsto \gamma(a,x)$ factors through $\theta$ for all $a \in \cA^n$. Since $\cA^n$ is a normal subgroup of $\Gammatil \ltimes \cA^n$, another application of Lemmas \ref{lem.dense-not-induced} and \ref{lem.further-trivialize} implies that $x \mapsto \gamma(g,x)$ factors through $\theta$ for all $g \in \Gammatil \ltimes \cA^n$.

We thus find a Borel map $\gamma_1 : (\Gammatil \ltimes \cA^n) \times \Gtil \to \Lambda$ such that for all $g \in \Gammatil \ltimes \cA^n$ and a.e.\ $x \in \Gtil \ltimes \R^n$, we have $\gamma(g,x) = \gamma_1(g,\theta(x))$.
From now on, we denote by $g$ the elements of $\Gammatil$ and we denote by $a$ the elements of $\cA^n$.
The restriction of $\gamma_1$ to $\Gammatil \times \Gtil$ is a $1$-cocycle for the translation action $\Gammatil \actson \Gtil$. Above, we have proven that this action is cocycle superrigid with countable target groups. Choose a Borel map $\vphi : \Gtil \to \Lambda$ and a group homomorphism $\delta : \Gammatil \to \Lambda$ such that $\gamma_1(g,y) = \vphi(gy)^{-1} \, \delta(g) \, \vphi(y)$ for all $g \in \Gammatil$ and a.e.\ $y \in \Gtil$.
Replacing $\gamma$ with the cohomologous $1$-cocycle $(g,x) \mapsto \vphi(\theta(gx)) \, \gamma(g,x) \, \vphi(\theta(x))^{-1}$, we may thus assume that $\gamma_1(g,y) = \delta(g)$ for all $g \in \Gammatil$ and a.e.\ $y \in \Gtil$. We denote $\delta_y(a) = \gamma_1(a,y)$ and note that $\delta_y : \cA^n \to \Lambda$ is a measurable family of group homomorphisms. To conclude the proof of the theorem, we have to show that for all $a \in \cA^n$, the map $y \mapsto \delta_y(a)$ is essentially constant.

When $g \in \Gammatil$, we have that $\pi(g) \in \SL(n,\cA^n)$ so that $\pi(g)(a) \in \cA^n$. The group law in $\Gammatil \ltimes \cA^n$ can then be expressed by $g \, a = \pi(g)(a) \, g$ for all $g \in \Gammatil$, $a \in \cA^n$. Applying the $1$-cocycle relation for $\gamma$, we conclude that
\begin{equation}\label{eq.important}
\delta_{gy}(\pi(g)a) = \delta(g) \, \delta_y(a) \, \delta(g)^{-1} \quad\text{for all $g \in \Gammatil$, $a \in \cA^n$ and a.e.\ $y \in \Gtil$.}
\end{equation}
Fix $i \in [n]$. Using the notation introduced above, denote by $L_i(\cA)$ the subgroup of $\Gamma$ generated by $H_i(\cA)$ and $R_i(\cA)$. Similarly define $L_i(\R)$. Then, $L_i(\cA)$ is dense in $L_i(\R)$ and $L_i(\R)$ consists of the matrices $A \in \SL(n,\R)$ satisfying $A(e_i) = e_i$, where $e_i$ is the $i$'th standard basis vector. Since the inclusion $\SL(n-1,\R) \cong H_i(\R) < \SL(n,\R)$ induces a surjective homomorphism between the fundamental groups, we get that $\pi^{-1}(L_i(\R))$ is a connected subgroup of $\Gtil$.

For every $a \in \cA$, denote $e_i(a) = a e_i \in \cA^n$. It follows from \eqref{eq.important} that for all $g \in \pi^{-1}(L_i(\cA))$, we have $\delta_{gy}(e_i(a)) = \delta(g) \, \delta_y(e_i(a)) \, \delta(g)^{-1}$. By Lemma \ref{lem.dense-not-induced}, for every fixed $i \in [n]$ and $a \in \cA$, the map $y \mapsto \delta_y(e_i(a))$ is invariant under left translation by the connected group $\pi^{-1}(L_i(\R))$, and thus of the form $y \mapsto \zeta(\pi(y)^{-1}(e_i))$ for some Borel map $\zeta : \R^n \to \Lambda$. This means that we find measurable families $(\rho_{i,z})_{z \in \R^n}$ of group homomorphisms $\rho_{i,z} : \cA \to \Lambda$ such that $\delta_y(e_i(a)) = \rho_{i,\pi(y)^{-1}(e_i)}(a)$ for all $i \in [n]$, $a \in \cA$ and a.e.\ $y \in \Gtil$.

Take $i \neq j$. We now apply \eqref{eq.important} for $g \in \Gammatil$ with  $\pi(g) = e_{ij}(-1)$ and $e_j(a) \in \cA^n$. We conclude that
$$\delta_{gy}(e_j(a)) \, \delta_{gy}(e_i(a))^{-1} = \delta_{gy}(e_j(a) - e_i(a)) = \delta_{gy}(\pi(g)e_j(a)) = \delta(g) \, \delta_y(e_j(a)) \, \delta(g)^{-1} \; .$$
Since $\pi(g)^{-1}(e_i) = e_i$ and $\pi(g)^{-1}(e_j) = e_i + e_j$, we find that
$$\rho_{j,\pi(y)^{-1}(e_i) + \pi(y)^{-1}(e_j)}(a) \, \rho_{i,\pi(y)^{-1}(e_i)}(a)^{-1} = \delta(g) \, \rho_{j,\pi(y)^{-1}(e_j)}(a) \, \delta(g)^{-1}$$
for all $a \in \cA$ and a.e.\ $y \in \Gtil$. Since $n \geq 3$, the map $\Gtil \to \R^n \times \R^n : y \mapsto (\pi(y)^{-1}(e_i),\pi(y)^{-1}(e_j))$ is a nonsingular factor map. We thus conclude that
\begin{equation}\label{eq.again-important}
\rho_{j,u+v}(a) \, \rho_{i,u}(a)^{-1} = \delta(g) \, \rho_{j,v}(a) \, \delta(g)^{-1} \quad\text{for all $a \in \cA$ and a.e.\ $(u,v) \in \R^n \times \R^n$.}
\end{equation}
Denote by $\cP$ the Polish group of Borel maps from $\R^n$ to $\Lambda$, where two such maps are identified if they are equal a.e., where the topology is given by convergence in measure and where the group law is defined pointwise. For every $v \in \R^n$ and $F \in \cP$, define $F_v \in \cP$ by $F_v(u) = F(u+v)$. Then, the map $\R^n \to \cP : v \mapsto F_v$ is continuous. Fix $a \in \cA$ and define $F,G \in \cP$ by $F(u) = \rho_{j,u}(a)$ and $G(u) = \rho_{i,u}(a)^{-1}$. Then, the map $\R^n \to \cP : v \mapsto F_v \, G$ is continuous. By \eqref{eq.again-important}, this map takes values in the discrete subgroup $\Lambda < \cP$ of constant functions. Since $\R^n$ is connected, it follows that $v \mapsto F_v$ is essentially constant. That means that we find group homomorphisms $\rho_j : \cA \to \Lambda$ such that $\rho_{j,u}(a) = \rho_j(a)$ for a.e.\ $y \in \R^n$.

We conclude that for all $j \in [n]$ and every $a \in \cA$, the map $y \mapsto \delta_y(e_j(a))$ is essentially constant. So, also $y \mapsto \delta_y(a)$ is essentially constant for every $a \in \cA^n$. This concludes the proof of the theorem.
\end{proof}

The first part of Theorem \ref{thm.main-no-go-theorem} can now be immediately deduced from Theorem \ref{thm.essential-cocycle-superrigidity}. We state and prove the following more general version.

\begin{corollary}\label{cor.an-example-of-OE-superrigid-v1}
Let $\cA \subset \R$ be any countable subring containing an algebraic number that does not belong to $\Z$. Let $\cF \subset \cA^*$ be a subgroup and $n \geq 3$ an integer. Define $\Gamma < \GL(n,\cA)$ as the group of matrices with $\det A \in \cF$. Consider the action of $G = \Gamma \ltimes \cA^n$ on $X=\R^n$ by $(A,a) \cdot x = A(a+x)$.

The action $G \actson X$ is essentially free, ergodic, nonsingular, simple and cocycle superrigid with countable targets. So the action is also OE superrigid (v1).

Denote by $T$ the closure of $\{ |a| \mid a \in \cF \}$ in $\R^*_+$. If $T = \R^*_+$, the action is of type III$_1$. If $T = \lambda^\Z$, the action is of type III$_\lambda$. If $T = \{1\}$, the action is of type II$_\infty$.
\end{corollary}

\begin{proof}
Since $\cA^n$ is dense in $\R^n$, by Lemma \ref{lem.dense-not-induced}, the action $\cA^n \actson \R^n$ is ergodic and not induced. A fortiori, $G \actson \R^n$ is ergodic and not induced. Assume that $\Sigma \lhd G$ is a normal subgroup whose action on $\R^n$ admits a fundamental domain. Write $\Sigma_0 = \Sigma \cap \cA^n$. Let $(A,a) \in \Sigma$. If $A \neq 1$, then for all $b \in \cA^n$,
$$(1,(1-A)b) = (1,b) (A,a) (1,b)^{-1} (A,a)^{-1} \in \Sigma$$
so that $\Sigma_0 \neq \{0\}$. Since $\Sigma_0$ is globally invariant under $\SL(n,\cA)$, the closure of $\Sigma_0$ in $\R^n$ is a nontrivial closed subgroup of $\R^n$ that is globally invariant under $\SL(n,\R)$. So, $\Sigma_0$ is dense in $\R^n$. It follows that $\Sigma_0 \actson \R^n$ is ergodic, contradicting the assumption that $\Sigma_0 \actson \R^n$ admits a fundamental domain. So, $A = 1$ and we have proven that $\Sigma \subset \cA^n$. If $\Sigma \neq \{0\}$, we again find that $\Sigma$ is dense in $\R^n$. So, $\Sigma$ is trivial and we have proven that $G \actson \R^n$ is a simple action.

Write $G_0 = \SL(n,\cA) \ltimes \cA^n$. By Theorem \ref{thm.essential-cocycle-superrigidity}, the dense subgroup $G_0$ of $\SL(n,\R) \ltimes \R^n$ is essentially cocycle superrigid with countable targets. By \cite[Proposition 3.3]{DV21}, the action $G_0 \actson \R^n$ is cocycle superrigid with countable targets. As mentioned above, this action is ergodic and not induced. Since $G_0$ is a normal subgroup of $G$, it follows from Lemma \ref{lem.further-trivialize} that also $G \actson \R^n$ is cocycle superrigid with countable targets. By Proposition \ref{prop.OE-v1-versus-cocycle}, the action is also OE superrigid (v1).

The Maharam extension of $G \actson \R^n$ can be identified with the action $G \actson \R^n \times \R$ given by
$$(A,a) \cdot (x,s) = (A(a+x),\log |\det A| + s) \; .$$
Since the translation action $\cA^n \actson \R^n$ is ergodic, it follows that the $G$-invariant functions on $\R^n \times \R$ are precisely the functions on $\R$ that are invariant under translation by all $\log |\det A|$, $A \in \GL(n,\cA)$, $\det A \in \cF$. So, these are the functions on $\R$ that are invariant under translation by $\{\log |a| \mid a \in \cF\}$, so that the type of $G \actson \R^n$ is as described in the corollary.
\end{proof}

\begin{proposition}\label{prop.no-go-III-0}
Let $G \actson (X,\mu)$ be a free, ergodic, nonsingular action of type III$_0$. If $G \actson X$ is not induced, then $G \actson X$ is not cocycle superrigid with countable targets. If $G \actson X$ is simple, then $G \actson (X,\mu)$ is not OE superrigid (v1).
\end{proposition}
\begin{proof}
Let $G \actson (X,\mu)$ be a free, ergodic, nonsingular action of type III$_0$. Assume that every $1$-cocycle $\Om : G \times X \to \Z$ is cohomologous to a group homomorphism. We prove that $G \actson X$ must be an induced action. By Proposition \ref{prop.OE-v1-versus-cocycle}, this suffices to prove the proposition.

Combining \cite[Theorem 2.7 and Remark 2.9]{Sch79} and \cite[Theorem 2.1]{JS85}, we find a free, ergodic, pmp action $\Z \actson (Y,\eta)$ and a Borel map $\pi : X \to Y$ such that $\pi_*(\mu) \sim \eta$ and $\pi(G \cdot x) \subset \Z \cdot \pi(x)$ for a.e.\ $x \in X$. Since all free, ergodic, pmp actions of $\Z$ are orbit equivalent, we may assume that the action $\Z \actson (Y,\eta)$ is the profinite $\Z \actson \Z_2$, viewed as the inverse limit of $\Z \actson \Z / 2^k \Z$, $k \in \N$.

Define the $1$-cocycle $\Om : G \times X \to \Z$ such that $\pi(g \cdot x) = \Om(g,x) \cdot \pi(x)$ for all $g \in G$ and a.e.\ $x \in X$. By our assumption, we find a Borel map $\vphi : X \to \Z$ and a group homomorphism $\delta : G \to \Z$ such that $\Om(g,x) = -\vphi(g \cdot x) + \delta(g) + \vphi(x)$ for all $g \in G$ and a.e.\ $x \in X$. Define the Borel map $\pi_1 : X \to Y : \pi_1(x) = \vphi(x) \cdot \pi(x)$. By construction, $\pi_1(g \cdot x) = \delta(g) \cdot \pi_1(x)$.

Note that $\pi_1$ is not essentially constant, since otherwise $\pi(x)$ takes values in a countable set for a.e.\ $x \in X$, contradicting $\pi_*(\mu) \sim \eta$. We can then choose $k \in \N$ large enough such that, denoting by $\psi : \Z_2 \to \Z / 2^k \Z$ the canonical quotient map, $\theta = \psi \circ \pi$ is not essentially constant. Since $\theta(g \cdot x) = \delta(g) + \theta(x)$ for all $g \in G$ and a.e.\ $x \in X$, it follows that the action $G \actson (X,\mu)$ is induced in a nontrivial way.
\end{proof}

\begin{proof}[{Proof of Theorem \ref{thm.main-no-go-theorem}}]
This now follows immediately from Corollary \ref{cor.an-example-of-OE-superrigid-v1} and Proposition \ref{prop.no-go-III-0}.
\end{proof}

\section{\boldmath Cocycle and OE superrigidity for actions of type III$_0$}

Let $G \actson (X,\mu)$ be a nonsingular action, with logarithm of the Radon-Nikodym $1$-cocycle $\om : G \times X \to \R$. Consider the Maharam extension $G \actson X \times \R$, together with the commuting measure scaling action $\R \actson X \times \R$ (see \eqref{eq.maharam} and \eqref{eq.commuting-scaling}). Denote by $\R \actson (Y,\eta)$ the associated flow. As explained at the start of Section \ref{sec.prelim}, we have a strictly $G$-invariant and $\R$-equivariant Borel map $\pi : X \times \R \to Y$ and we define $\psi : X \to Y$ by \eqref{eq.choice-psi}.
%


Given a subgroup $\Lambda_0 < \Lambda$, we denote by $C_\Lambda(\Lambda_0) = \{g \in \Lambda \mid gh = hg \;\text{for all $h \in \Lambda_0$}\}$ the centralizer of $\Lambda_0$ inside $\Lambda$.

\begin{theorem}\label{thm.key-cocycle-superrigid}
Let $G \actson (X,\mu)$ be a nonsingular action of a countable group $G$ on a standard probability space $(X,\mu)$. Let $\om : G \times X \to \R$ be the logarithm of the Radon-Nikodym cocycle and let $G \actson (\Xtil,\mutil)$ be the Maharam extension, with ergodic decomposition $(\Xtil_y,\mutil_y)_{y \in Y}$, associated flow $\R \actson (Y,\eta)$ and Borel map $\psi : X \to Y$ as in \eqref{eq.choice-psi}.

Let $\Lambda$ be a countable group. Assume that $G$ is finitely generated and that for $\eta$-a.e.\ $y \in Y$, the ergodic action $G \actson (\Xtil_y,\mutil_y)$ is not induced and cocycle superrigid with target $\Lambda$. Then for any $1$-cocycle $\Om : G \times X \to \Lambda$, there exists a group homomorphism $\delta : G \to \Lambda$ and a strict $1$-cocycle $\gamma : \R \times Y \to C_\Lambda(\delta(G))$ such that $\Om$ is cohomologous with the $1$-cocycle
    $$G \times X \to \Lambda : (g,x) \mapsto \delta(g) \, \gamma(-\om(g,x),\psi(x)) \; .$$
%
\end{theorem}
\begin{proof}
%
%
Define the $1$-cocycle $\Omtil : (G \times \R) \times \Xtil \to \Lambda : \Omtil((g,t),(x,s)) = \Om(g,x)$. First restrict $\Omtil$ to a $1$-cocycle for the action $G \actson \Xtil$ with ergodic decomposition given by $\pi : \Xtil \to Y$. As explained in detail in \cite{FMW04}, we may consider $\Omtil$ as a measurable family $(\Omtil_y)_{y \in Y}$ of $1$-cocycles for the measurable family of actions $G \actson (\Xtil_y,\mutil_y)$. By assumption, $\eta$-a.e.\ $\Omtil_y$ is cohomologous to a group homomorphism $\delta_y : G \to \Lambda$.

By \cite[Corollary 3.11]{FMW04}, we find a Borel family of group homomorphisms $\delta_y : G \to \Lambda$, indexed by $y \in Y$, and a Borel map $\vphi : \Xtil \to \Lambda$ such that for all $g \in G$, we have that
$$\Omtil((g,0),(x,s)) = \vphi((g,0) \cdot (x,s))^{-1} \, \delta_{\pi(x,s)}(g) \, \vphi(x,s) \quad\text{for $\mutil$-a.e.\ $(x,s) \in \Xtil$.}$$
Define the $1$-cocycle $\Psi : (G \times \R) \times \Xtil \to \Lambda$ by
$$\Psi((g,t),(x,s)) = \vphi((g,t) \cdot (x,s)) \, \Omtil((g,t),(x,s)) \, \vphi(x,s)^{-1} \; .$$
By construction, $\Psi \sim \Omtil$ as $1$-cocycles for $G \times \R \actson \Xtil$ and for all $g \in G$, $\Psi((g,0),(x,s)) = \delta_{\pi(x,s)}(g)$ for $\mutil$-a.e.\ $(x,s) \in \Xtil$.
%
Define $\zeta_t(x,s) = \Psi((e,t),(x,s))$. From the $1$-cocycle relation for $\Psi$ applied to $(g,0)(e,t) = (g,t) = (e,t)(g,0)$, it follows that for all $t \in \R$, $g \in G$, we have
\begin{equation}\label{eq.one}
\zeta_t(g \cdot (x,s)) = \delta_{t \cdot \pi(x,s)}(g) \, \zeta_t(x,s) \, \delta_{\pi(x,s)}(g)^{-1} \quad\text{for $\mutil$-a.e.\ $(x,s) \in \Xtil$.}
\end{equation}
Fix $t \in \R$. Then \eqref{eq.one} is saying that for $\eta$-a.e.\ $y \in Y$, $\zeta_t$ is a $G$-equivariant Borel map from $(\Xtil_y,\mutil_y)$ to the countable set $\Lambda$ on which $G$ is acting by $g \cdot \lambda = \delta_{t \cdot y}(g) \lambda \delta_y(g)^{-1}$. Since we assumed that for $\eta$-a.e.\ $y \in Y$, the action $G \actson (\Xtil_y,\mutil_y)$ is not induced, it follows that for $\eta$-a.e.\ $y \in Y$, the map $\zeta_t$ is $\mutil_y$-a.e.\ constant on $\Xtil_y$.

We thus find a Borel map $\gamma_0 : \R \times Y \to \Lambda$ such that for all $t \in \R$, we have that $\zeta_t(x,s) = \gamma_0(t,\pi(x,s))$ for $\mutil$-a.e.\ $(x,s) \in \Xtil$. Since $\Omtil((e,t),(x,s)) = e$, we have $\zeta_t(x,s) = \vphi(x,t+s) \, \vphi(x,s)^{-1}$ for $\mutil$-a.e.\ $(x,s) \in \Xtil$. Therefore, for every $t \in \R$,
\begin{equation}\label{eq.two}
\gamma_0(t,\pi(x,s)) = \vphi(x,t+s) \, \vphi(x,s)^{-1} \quad\text{for $\mutil$-a.e.\ $(x,s) \in \Xtil$.}
\end{equation}
So, $\gamma_0 : \R \times Y \to \Lambda$ is a $1$-cocycle. Then \eqref{eq.one} is saying that for all $t \in \R$, $g \in G$, we have
\begin{equation}\label{eq.three}
\delta_{t \cdot y}(g) = \gamma_0(t,y) \, \delta_y(g) \, \gamma_0(t,y)^{-1} \quad\text{for $\eta$-a.e.\ $y \in Y$.}
\end{equation}

Since $G$ is finitely generated, the set of group homomorphisms from $G$ to $\Lambda$ is countable. We thus find a group homomorphism $\delta : G \to \Lambda$ such that $\delta_y = \delta$ for all $y$ in a nonnegligible Borel subset $\cU \subset Y$. Combining \eqref{eq.three} with the ergodicity of $\R \actson (Y,\eta)$, it follows that $\delta_y$ is conjugate to $\delta$ for $\eta$-a.e.\ $y \in Y$. We then find a Borel map $\rho : Y \to \Lambda$ so that $\delta_y(g) = \rho(y)^{-1} \, \delta(g) \, \rho(y)$ for $\eta$-a.e.\ $y \in Y$ and all $g \in G$. Replacing $\vphi(x,s)$ by $\rho(\pi(x,s)) \vphi(x,s)$, we may assume that $\delta_y = \delta$ for $\eta$-a.e.\ $y \in Y$.

The $1$-cocycle $\Psi : (G \times \R) \times \Xtil \to \Lambda$ thus has the property that for all $g \in G$ and $t \in \R$, we have that
$$\Psi((g,0),(x,s)) = \delta(g) \quad\text{and}\quad \Psi((e,t),(x,s)) = \gamma_0(t,\pi(x,s))$$
for $\mutil$-a.e.\ $(x,s) \in \Xtil$. The cocycle identity for $\Psi$ then forces $\gamma_0$ to take values a.e.\ in the centralizer $C_\Lambda(\delta(G))$.

Choose a strict $1$-cocycle $\gamma : \R \times Y \to C_\Lambda(\delta(G))$ such that for every $t \in \R$, we have that $\gamma(t,y) = \gamma_0(t,y)$ for $\eta$-a.e.\ $y \in Y$.
Define the Borel map $\theta : \Xtil \to C_\Lambda(\delta(G)) : \theta(x,s) = \gamma(s,\psi(x))$. Consider the cohomologous $1$-cocycle $\Psi_1 \sim \Psi$ defined by
$$\Psi_1((g,t),(x,s)) = \theta((g,t)\cdot(x,s))^{-1} \, \Psi((g,t),(x,s)) \, \theta(x,s) \; .$$
Since $\gamma$ is a strict $1$-cocycle, we find for every $t \in \R$ and $g \in G$ that $\Psi_1((e,t),(x,s)) = e$ and $\Psi_1((g,0),(x,s)) = \delta(g) \, \gamma(-\om(g,x),\psi(x))$ for $\mutil$-a.e.\ $(x,s) \in \Xtil$. So, defining the $1$-cocycle
$$\Psi_0 : G \times X \to \Lambda : \Psi_0(g,x) = \delta(g) \, \gamma(-\om(g,x),\psi(x)) \; ,$$
we have proven that the $1$-cocycles $\Om$ and $\Psi_0$ are cohomologous when viewed as $1$-cocycles for $G \times \R \actson \Xtil$. The $\R$-invariance of both $1$-cocycles forces the Borel function implementing the cohomology $\Om \sim \Psi_0$ to be essentially $\R$-invariant as well. We have thus proven that $\Om \sim \Psi_0$.
\end{proof}

\begin{remark}\label{rem.reformulation}
The conclusion of Theorem \ref{thm.key-cocycle-superrigid} can also be formulated in the following way. Consider the Maharam extension, together with its measure scaling action, $G \times \R \actson \Xtil = X \times \R$. The conclusion of Theorem \ref{thm.key-cocycle-superrigid} says that the $1$-cocycles for $G \times \R \actson X \times \R$ given by
$$((g,r),(x,t)) \mapsto \Om(g,x) \quad\text{and}\quad ((g,r),(x,t)) \mapsto \delta(g) \, \gamma(-\om(g,x),\psi(x)) \; ,$$
and which are both trivial on $\R$, are cohomologous. But using the map $(x,t) \mapsto \gamma(t,\psi(x))$, this second $1$-cocycle is also cohomologous to $((g,r),(x,t)) \mapsto \delta(g) \, \gamma(r,\pi(x,t))$, where $\pi : X \times \R \to Y$ is the ergodic decomposition of $G \actson X \times \R$.
\end{remark}

Theorem \ref{thm.key-cocycle-superrigid} applies in particular to the type III actions of the form \eqref{eq.our-actions}. This then leads to the following cocycle superrigidity result.

\begin{theorem}\label{thm.key-cocycle-superrigid-v2}
Let $G \actson (X,\mu)$ be a free, ergodic, nonsingular action of type III$_1$. Denote by $\om : G \times X \to \R$ the logarithm of the Radon-Nikodym cocycle. Let $\R \actson (Y,\eta)$ be an ergodic flow and consider
\begin{equation}\label{eq.our-action-repeated}
G \actson X \times Y : g \cdot (x,y) = (g \cdot x, \om(g,x) \cdot y)
\end{equation}
as in \eqref{eq.our-actions}.

If $G$ is finitely generated and if the Maharam extension of $G \actson (X,\mu)$ is not induced and cocycle superrigid with countable targets, then every $1$-cocycle $\Om : G \times X \times Y \to \Lambda$ for the action \eqref{eq.our-action-repeated} with values in a countable group $\Lambda$ is cohomologous with a $1$-cocycle of the form
$$G \times X \times Y \to \Lambda : (g,x,y) \mapsto \delta(g) \, \gamma(\om(g,x),y) \; ,$$
where $\delta : G \to \Lambda$ is a group homomorphism and $\gamma : \R \times Y \to C_\Lambda(\delta(G))$ is a $1$-cocycle.
\end{theorem}
\begin{proof}
As in the proof of Proposition \ref{prop.prescribed-flow}, we take the unique action $\R^2 \actson (Z,\zeta)$ given by Proposition \ref{prop.unique-char-adj}, associated with the ergodic flow $\R \actson (Y,\eta)$. Identify $Y = Z / (\{0\} \times \R)$ and denote by $\pi_1 : Z \to Y$ the corresponding factor map. Also write $\Yh = Z/(\R \times \{0\})$ and denote by $\pi_2 : Z \to \Yh$ the corresponding factor map. The Maharam extension of $G \actson X \times Y$ together with its measure scaling action of $\R$ is then given by
$$G \times \R \actson X \times Z : (g,t) \cdot (x,z) = (g\cdot x, (\om(g,x),t) \cdot z) \; .$$
The map $(x,z) \mapsto \pi_2(z)$ identifies the associated flow of the action $G \actson X \times Y$ with $\R \actson \Yh$. Identifying $Z = \R \times \Yh$, the ergodic decomposition of $G \actson X \times Z$ is a.e.\ given by the Maharam extension $G \actson X \times \R$ of the initial type III$_1$ action $G \actson X$. So, $G \actson X \times Y$ satisfies the assumptions of Theorem \ref{thm.key-cocycle-superrigid}.

Let $\Lambda$ be a countable group and $\Om : G \times X \times Y \to \Lambda$ a $1$-cocycle. Define the $1$-cocycle
$$\Omtil : (G \times \R) \times (X \times Z) \to \Lambda : \Omtil((g,t),(x,z)) = \Om(g,(x,\pi_1(z))) \; .$$
It follows from Theorem \ref{thm.key-cocycle-superrigid} and Remark \ref{rem.reformulation} that $\Omtil$ is cohomologous with $\Om_1$, where
$$\Om_1((g,t),(x,z)) = \delta(g) \, \gamma_1(t,\pi_2(z)) \; ,$$
with $\delta : G \to \Lambda$ a group homomorphism and $\gamma_1 : \R \times \Yh \to C_\Lambda(\delta(G))$ a $1$-cocycle. Define the $1$-cocycle
$$\gammatil_1 : \R^2 \times Z  \to C_\Lambda(\delta(G)) : \gammatil_1((r,t),z) = \gamma_1(t,\pi_2(z)) \; .$$
Since the action of $\{0\} \times \R$ on $Z$ is measure scaling, $\gammatil_1$ is cohomologous to a $1$-cocycle $\gamma_2$ of the form $\gamma_2((r,t),z) = \gamma(r,\pi_1(z))$, where $\gamma : \R \times Y \to C_\Lambda(\delta(G))$ is a $1$-cocycle. Choose a Borel map $\vphi : Z \to C_\Lambda(\delta(G))$ implementing this cohomology, so that for all $(r,t) \in \R^2$, we have
$$\gamma(r,\pi_1(z)) = \vphi((r,t) \cdot z) \, \gamma_1(t,\pi_2(z)) \, \vphi(z)^{-1} \quad\text{for a.e.\ $z \in Z$.}$$
Define the $1$-cocycle $\Om_2 : (G \times \R) \times (X \times Z) \to \Lambda$ by
$$\Om_2((g,t),(x,z)) = \vphi((\om(g,x),t) \cdot z) \, \Om_1((g,t),(x,z)) \, \vphi(z)^{-1} \; .$$
By construction, $\Omtil \sim \Om_2$ and
$$\Om_2((g,t),(x,z)) = \delta(g) \, \gamma(\om(g,x),\pi_1(z)) \; .$$
Since both $\Omtil$ and $\Om_2$ are trivial on $\{0\} \times \R$, this means that $\Om$ is cohomologous with the $1$-cocycle
$$G \times X \times Y : (g,x,y) \mapsto \delta(g) \, \gamma(\om(g,x),y) \; .$$
\end{proof}

\begin{remark}
When $G \actson (X,\mu)$ is a free, ergodic, nonsingular action of type III, with Maharam extension $G \actson (\Xtil,\mutil)$ whose ergodic decomposition is denoted as $(\Xtil_y,\mutil_y)_{y \in Y}$, it follows from \cite[Theorem XII.1.1]{Tak03} that for a.e.\ $y \in Y$, the action $G \actson (\Xtil_y,\mutil_y)$ is of type II$_\infty$. To give examples where Theorem \ref{thm.key-cocycle-superrigid} applies, we thus need cocycle superrigidity for concrete actions of type II$_\infty$. For the specific actions appearing in Theorem \ref{thm.key-cocycle-superrigid-v2}, by construction, the actions $G \actson (\Xtil_y,\mutil_y)$ are a.e.\ the same.
\end{remark}

Both Theorem \ref{thm.key-cocycle-superrigid} and \ref{thm.key-cocycle-superrigid-v2} immediately lead to OE superrigidity results. We start with the following result.

\begin{corollary}\label{cor.key-OE-superrigid}
Let $G$ be a finitely generated group with trivial center and let $G \actson (X,\mu)$ be an essentially free, nonsingular, ergodic action. Let $\om : G \times X \to \R$ be the logarithm of the Radon-Nikodym cocycle and let $G \actson (\Xtil,\mutil)$ be the Maharam extension, with associated flow $\R \actson (Y,\eta)$ and Borel map $\psi : X \to Y$ as in \eqref{eq.choice-psi}. Assume that for $\eta$-a.e.\ $y \in Y$, the action $G \actson (\Xtil_y,\mutil_y)$ is simple and cocycle superrigid with countable target groups.

Then, $G \actson (X,\mu)$ satisfies the OE-superrigidity property (v2) defined in the introduction.

More precisely, any free nonsingular ergodic action that is stably orbit equivalent with $G \actson (X,\mu)$ is conjugate to an induction of an action of the form
\begin{equation}\label{eq.precise-action}
G \times \Lambda \actson X \times \Lambda : (g,a) \cdot (x,b) = (g \cdot x, \gamma(-\om(g,x),\psi(x))b a^{-1}) \; ,
\end{equation}
where $\Lambda$ is a countable group and $\gamma : \R \times Y \to \Lambda$ is a strict $1$-cocycle.
\end{corollary}


\begin{proof}
Let $\Gamma \actson (Z,\zeta)$ be a free, ergodic, nonsingular action and let $\Delta : \cU \subset X \to Z$ be a stable orbit equivalence between $G \actson X$ and $\Gamma \actson Z$. By ergodicity of $G \actson (X,\mu)$, we can choose a Borel map $\theta : X \to G$ such that $\theta(x) = e$ for all $x \in \cU$ and $\theta(x) \cdot x \in \cU$ for a.e.\ $x \in X$. Define $\Delta_0 : X \to Z : \Delta_0(x) = \Delta(\theta(x) \cdot x)$. We then define the Zimmer $1$-cocycle $\Om : G \times X \to \Gamma$ such that $\Delta_0(g \cdot x) = \Om(g,x) \cdot \Delta_0(x)$ for all $g \in G$ and a.e.\ $x \in X$.

To translate the cocycle superrigidity theorem \ref{thm.key-cocycle-superrigid} to an OE superrigidity theorem, we use the connection with measure equivalence as developed in \cite[Section 3]{Fur98} (see also \cite[Lemma 2.2]{DV21} for a result that exactly suits our purposes). Define the action
\begin{equation}\label{eq.first-action}
G \actson X \times \Gamma : g \cdot (x,b) = (g \cdot x, \Om(g,x) b) \; ,
\end{equation}
which commutes with the right translation action by $\Gamma$ in the second variable. By the results cited above, the action $G \actson X \times \Gamma$ admits a fundamental domain and there is a natural isomorphism of $\Gamma$-actions $\al : G \backslash (X \times \Gamma) \to Z$ with the property that $\Delta(x) \in \Gamma \cdot \al(x,e)$ for a.e.\ $x \in \cU$.

By Theorem \ref{thm.key-cocycle-superrigid}, we find a group homomorphism $\delta : G \to \Gamma$ and a $1$-cocycle $\gamma : \R \times Y \to C_\Gamma(\delta(G))$ such that $\Om$ is cohomologous with the $1$-cocycle
$$\Om_1 : G \times X \to \Gamma : \Om_1(g,x) = \delta(g) \, \gamma(-\om(g,x),\psi(x)) \; .$$
Let $\vphi : X \to \Gamma$ be a Borel map such that $\Om_1(g,x) = \vphi(g \cdot x) \, \Om(g,x) \, \vphi(x)^{-1}$. The map $(x,b) \mapsto (x,\vphi(x)b)$ implements an isomorphism between the action $G \actson X \times \Gamma$ in \eqref{eq.first-action} and the action
\begin{equation}\label{eq.second-action}
G \actson X \times \Gamma : g \cdot (x,b) = (g \cdot x, \Om_1(g,x) b) \; .
\end{equation}
Moreover, the action and the isomorphism commute with the $\Gamma$-action. We thus still find an isomorphism of $\Gamma$-actions $\al_1 : G \backslash (X \times \Gamma) \to Z$ with the property that $\Delta(x) \in \Gamma \cdot \al_1(x,e)$ for a.e.\ $x \in \cU$.

The $1$-cocycle $\Omtil_1 : G \times \Xtil : \Omtil_1(g,(x,s)) = \Om_1(g,x)$ for the Maharam extension $G \actson \Xtil$ is, by construction, cohomologous with the $1$-cocycle $(g,(x,s)) \mapsto \delta(g)$. Since the action $G \actson X \times \Gamma$ admits a fundamental domain, a fortiori, the same holds for the action $G \actson \Xtil \times \Gamma : g \cdot (x,s,b) = (g \cdot (x,s), \delta(g)b)$, and thus for the action $\Ker \delta \actson \Xtil$. Since we assumed that a.e.\ action $G \actson \Xtil_y$ is simple, the normal subgroup $\Ker \delta$ must be trivial. So, $\delta : G \to \Gamma$ is faithful. Define $\Lambda = C_\Gamma(\delta(G))$. Since $G$ has trivial center, $\delta(G) \cap \Lambda = \{e\}$. We have thus found a subgroup $\delta(G) \times \Lambda < \Gamma$.

Since $\Om_1$ takes values in $\delta(G) \times \Lambda$, the action $\Gamma \actson G \backslash (X \times \Gamma)$ is induced from $\delta(G) \times \Lambda \actson G \backslash (X \times \delta(G) \times \Lambda)$. Under the natural identification  $G \backslash (X \times \delta(G) \times \Lambda) = X \times \Lambda$ and the isomorphism $\delta \times \id : G \times \Lambda \to \delta(G) \times \Lambda$, this last action is precisely the action given by \eqref{eq.precise-action}.
\end{proof}


In exactly the same way as Theorem \ref{thm.key-cocycle-superrigid-v2} is deduced from Theorem \ref{thm.key-cocycle-superrigid}, we can deduce the following result from Corollary \ref{cor.key-OE-superrigid}. We thus omit the proof. Note that Theorem \ref{thm.main-OE-superrigid} is contained in the following corollary.

\begin{corollary}\label{cor.OE-superrigid-prescribed-flow}
Let $G \actson (X,\mu)$ be a free, ergodic, nonsingular action of type III$_1$. Assume that $G$ is finitely generated and has trivial center. Assume that the Maharam extension of $G \actson (X,\mu)$ is simple and cocycle superrigid with countable targets. Denote by $\om : G \times X \to \R$ the logarithm of the Radon-Nikodym cocycle. Let $\R \actson (Y,\eta)$ be an ergodic flow and consider
\begin{equation}\label{eq.our-actions-again-repeated}
G \actson X \times Y : g \cdot (x,y) = (g \cdot x, \om(g,x) \cdot y)
\end{equation}
as in \eqref{eq.our-actions}.

Then, $G \actson (X \times Y,\mu \times \eta)$ satisfies the OE-superrigidity property (v2) defined in the introduction.

More precisely, any free nonsingular ergodic action that is stably orbit equivalent with the action \eqref{eq.our-actions-again-repeated} is conjugate to an induction of an action of the form
\begin{equation}\label{eq.my-new-action}
G \times \Lambda \actson X \times Y \times \Lambda : (g,a) \cdot (x,y,b) = (g \cdot x, \om(g,x) \cdot y, \gamma(\om(g,x),y) b a^{-1}) \; ,
\end{equation}
where $\Lambda$ is a countable group and $\gamma : \R \times Y \to \Lambda$ is a $1$-cocycle.
\end{corollary}

There are several concrete group actions that satisfy the assumptions of Corollary \ref{cor.OE-superrigid-prescribed-flow}. We start with the following example of \cite{PV08}.

\begin{example}\label{ex.the-ex-from-PV}
Whenever $n \geq 5$ is an odd integer and $G < \SL(n,\R)$ is a lattice, the action $G \actson \R^n/\R^*_+$ satisfies the hypotheses of Corollary \ref{cor.OE-superrigid-prescribed-flow}. Indeed, by \cite[Theorem 1.3]{PV08}, the Maharam extension $G \actson \R^n$ is cocycle superrigid. Moreover, these groups $G$ have property~(T), so that they are finitely generated. Since $n$ is odd, $G$ has trivial center. By \cite[Lemmas 5.6 and 6.1]{PV08}, the action $G \actson \R^n$ is doubly ergodic and not induced. By Margulis' normal subgroup theorem, a normal subgroup of $G$ is either trivial (because $n$ is odd) or of finite index, and thus acting ergodically on $\R^n$. So, the action $G \actson \R^n$ is simple.
\end{example}

In the following theorem, we prove that all the assumptions of Corollary \ref{cor.OE-superrigid-prescribed-flow} are satisfied for the action $G \actson \R^n / \R^*_+$ when $G$ ranges over a broad family of dense subgroups $\SL(n,\R)$ and $n \geq 3$ is an odd integer. This then leads to the proof of point~1 in Corollary \ref{cor.main-example} (see Corollary \ref{cor.main-example-with-out} below).

\begin{theorem}\label{thm.rings-cocycle-superrigid}
Let $\cA \subset \R$ be any countable subring containing an algebraic number that does not belong to $\Z$. Let $n \geq 3$ be an integer and $E(n,\cA) < G < \SL(n,\cA)$. Then the linear action $G \actson \R^n$ is essentially free, ergodic, nonsingular, not induced and cocycle superrigid with countable targets. If $n$ is odd, $G$ has trivial center and the action is simple. If $\cA$ is finitely generated as a ring, then $E(n,\cA)$ is a finitely generated group.
\end{theorem}
\begin{proof}
Since $G < \SL(n,\R)$ is dense, it follows from Lemma \ref{lem.dense-not-induced} that the action $G \actson \R^n$ is ergodic and not induced. By density of $G < \SL(n,\R)$, the center of $G$ belongs to the center of $\SL(n,\R)$, which is trivial if $n$ is odd. If $n$ is odd and $\Sigma \lhd G$ is a normal subgroup whose action on $\R^n$ admits a fundamental domain, the closure $\overline{\Sigma}$ of $\Sigma$ in $\SL(n,\R)$ is, by density of $G$, a normal subgroup of $\SL(n,\R)$. Since $n$ is odd, it follows that either $\overline{\Sigma} = \{1\}$, or $\overline{\Sigma} = \SL(n,\R)$. In the second case, $\Sigma$ is a dense subgroup of $\SL(n,\R)$, so that $\Sigma \actson \R^n$ is ergodic. It thus follows that $\Sigma = \{1\}$.

By Theorem \ref{thm.essential-cocycle-superrigidity}, $G < \SL(n,\R)$ is essentially cocycle superrigid with countable targets. By \cite[Proposition 3.3]{DV21}, the action $G \actson \R^n$ is cocycle superrigid with countable targets.

By \cite[Proposition 4.3.11]{HOM89}, the group $E(n,\cA)$ is finitely generated when $n \geq 3$ and $\cA$ is finitely generated as a ring.
\end{proof}


\section{Conjugacy and classification results}

In Corollary \ref{cor.key-OE-superrigid}, we proved that free, nonsingular, ergodic actions $G \actson (X,\mu)$ with a sufficiently rigid Maharam extension $G \actson (\Xtil,\mutil)$ satisfy the OE-superrigidity property (v2) and we described all stably orbit equivalent actions. We now prove the following complete classification up to conjugacy of this class of stably orbit equivalent actions.

For the formulation of the following proposition, note that every conjugacy of actions (and actually every stable orbit equivalence) gives rise to a canonical associated isomorphism between the associated flows.

\begin{proposition}\label{prop.complete-classification}
Let $G \actson (X,\mu)$ be an essentially free, nonsingular, ergodic action. Make the same assumptions as in Corollary \ref{cor.key-OE-superrigid}. Whenever $\gamma : \R \times Y \to \Lambda$ is a strict $1$-cocycle with values in a countable group $\Lambda$ and whenever $\sigma : G \times \Lambda \to \Gamma$ is a faithful group homomorphism to a countable group $\Gamma$, we denote by $\be(\gamma,\sigma)$ the $\Gamma$-action defined as the induction of the action
\begin{equation}\label{eq.precise-action-repeated}
\beta(\gamma) : G \times \Lambda \actson X \times \Lambda : (g,a) \cdot (x,b) = (g \cdot x, \gamma(-\om(g,x),\psi(x)) b a^{-1})
\end{equation}
along the embedding $\si : G \times \Lambda \to \Gamma$.
\begin{enumlist}
\item An essentially free, nonsingular, ergodic action is stably orbit equivalent with $G \actson (X,\mu)$ if and only if it is conjugate to $\be(\gamma,\si)$ for some $\gamma,\si$ as above.
\item The actions $\be(\gamma,\si)$ and $\be(\gamma',\si')$ are conjugate if and only if there exist subgroups $\Lambda_0 < \Lambda$, $\Lambda'_0 < \Lambda'$, an automorphism $\delta \in \Aut(G)$, group isomorphisms $\rho : \Lambda_0 \to \Lambda'_0$ and $\al : \Gamma \to \Gamma'$, and a $\delta$-conjugacy $\Delta : X \to X'$ with associated isomorphism $\Delta_0 : Y \to Y'$ of flows such that
    \begin{itemlist}
    \item $\gamma,\gamma'$ are cohomologous to strict $1$-cocycles $\gamma_0,\gamma'_0$ that take values in $\Lambda_0,\Lambda'_0$~;
    \item $\al(\si(g,a)) = \si'(\delta(g),\rho(a))$ for all $g \in G$ and $a \in \Lambda_0$~;
    \item the $1$-cocycles $\rho \circ \gamma_0$ and $\gamma'_0 \circ (\id \times \Delta_0)$ are cohomologous as $1$-cocycles $\R \times Y \to \Lambda_0'$.
    \end{itemlist}
\end{enumlist}
\end{proposition}

Point~1 of Proposition \ref{prop.complete-classification} is just a repetition of Corollary \ref{cor.key-OE-superrigid}. We deduce point~2 of Proposition \ref{prop.complete-classification} from the following two results. First in Proposition \ref{prop.when-induced}, we describe when and how an action of the form $\beta(\gamma)$ in \eqref{eq.precise-action-repeated} is induced. Second in Proposition \ref{prop.when-conjugate}, we prove when two actions of the form $\beta(\gamma)$ are conjugate.

For our main family of group actions $G \actson X \times Y$ defined in \eqref{eq.our-actions}, it then remains to analyze when two such actions are conjugate. Under the appropriate assumption, we prove in Proposition \ref{prop.special-when-conjugate} that this happens if and only if the $G$-actions $G \actson X$ are conjugate and the flows $\R \actson Y$ are isomorphic. In particular, we find the outer automorphism groups of these type III$_0$ orbit equivalence relations.

Before proving Proposition \ref{prop.complete-classification}, we clarify the following subtle point. When $\gamma$ and $\gamma'$ are strict $1$-cocycles that are cohomologous, expressed by the a.e.\ equality $\gamma'(t,y) = \vphi(t \cdot y) \, \gamma(t,y) \, \vphi(y)^{-1}$, there is a natural isomorphism $\Delta_{\gamma',\gamma}$ between the actions $\be(\gamma)$ and $\be(\gamma')$. This follows immediately by observing that the Maharam extension of $\be(\gamma)$, together with its commuting $\R$-action, is isomorphic with
\begin{equation}\label{eq.Maharam-of-beta-gamma}
G \times \Lambda \times \R \actson X \times \R \times \Lambda : (g,a,t) \cdot (x,s,b) = (g \cdot x, \om(g,x) + t + s, \gamma(t,\pi(x,s))ba^{-1}) \; .
\end{equation}
Moreover, the map
\begin{equation}\label{eq.map-Phi}
\Phi : X \times \R \times \Lambda \to X \times \Lambda : \Phi(x,s,b) = (x,\gamma(s,\psi(x))^{-1} b)
\end{equation}
is $\R$-invariant and $(G \times \Lambda)$-equivariant.

Then the map $(x,s,b) \mapsto (x,s,\vphi(\pi(x,s))b)$ is a well defined isomorphism between the Maharam extensions. Taking the quotient by the action of $\R$, we find $\Delta_{\gamma',\gamma}$.

%
%
%
%

\begin{proposition}\label{prop.when-induced}
Let $G \actson (X,\mu)$ be an essentially free, nonsingular, ergodic action with Maharam extension $G \actson (\Xtil,\mutil)$ and associated flow $\R \actson (Y,\eta)$. Assume that for $\eta$-a.e.\ $y \in Y$, the action $G \actson (\Xtil_y,\mutil_y)$ is not induced. Let $\gamma : \R \times Y \to \Lambda$ be a strict $1$-cocycle with values in a countable group $\Lambda$.

The action $\beta(\gamma)$ in \eqref{eq.precise-action-repeated} is induced from a subgroup $\Gamma < G \times \Lambda$ acting on $Z \subset X \times \Lambda$ if and only if $\Gamma = G \times \Lambda_0$ for a subgroup $\Lambda_0 < \Lambda$, $\gamma$ is cohomologous with a strict $1$-cocycle $\gamma_0$ taking values in $\Lambda_0$ and $\Delta_{\gamma_0,\gamma}(Z) = X \times \Lambda_0$.
%
\end{proposition}


\begin{proof}
If $\gamma$ takes values in $\Lambda_0$, we have by construction that $\be(\gamma)$ is induced from $G \times \Lambda_0$ acting on $X \times \Lambda_0$. So we only prove the converse and assume that $\be(\gamma)$ is induced from $\Gamma \actson Z$. Write $I = (G \times \Lambda)/\Gamma$. Consider the Maharam extension of $\be(\gamma)$ given by \eqref{eq.Maharam-of-beta-gamma}. Since $\be(\gamma)$ is induced from $\Gamma \actson Z$, we find a Borel map $\theta : X \times \R \times \Lambda \to I$ that is $\R$-invariant and $(G \times \Lambda)$-equivariant.

Since for a.e.\ $y \in Y$, we have that $G \actson \Xtil_y$ is not induced, the $G$-equivariance of $\theta$ implies that $\theta(x,s,b) = \theta_1(\pi(x,s),b)$, where $\theta_1 : Y \times \Lambda \to I$. In particular, $\theta$ is a $G$-invariant map. It follows that $G \times \{e\} \subset \Gamma$, so that $\Gamma = G \times \Lambda_0$ for some subgroup $\Lambda_0 < \Lambda$. From now on, we identify $I = \Lambda / \Lambda_0$, on which $G$ acts trivially.

The $\Lambda$-equivariance of $\theta$ implies that also $\theta_1$ is $\Lambda$-equivariant and thus, of the form $\theta_1(y,b) = b^{-1} \theta_2(y) \Lambda_0$ for some Borel map $\theta_2 : Y \to \Lambda$. Expressing the $\R$-invariance of $\theta$ and thus, the invariance of $\theta_1$ under the action $t \cdot (y,s) = (t \cdot y, \gamma(t,y) s)$, we find that
$$b^{-1} \theta_2(y) \Lambda_0 = \theta_1(y,b) = \theta_1(t \cdot y, \gamma(t,y) b) = b^{-1} \gamma(t,y)^{-1} \theta_2(t \cdot y) \Lambda_0 \; .$$
This precisely means that the cohomologous $1$-cocycle $\gamma_0(t,y) = \theta_2(t \cdot y)^{-1} \gamma(t,y) \theta_2(y)$ takes values in $\Lambda_0$. Using the notation introduced before the proposition, this also means that $\Delta_{\gamma_0,\gamma}(Z) = X \times \Lambda_0$.
\end{proof}

\begin{proposition}\label{prop.when-conjugate}
For $i \in \{1,2\}$, let $G_i \actson (X_i,\mu_i)$ be essentially free, nonsingular, ergodic actions with Maharam extension $G_i \actson (\Xtil_i,\mutil_i)$ and associated flow $\R \actson (Y_i,\eta_i)$. Assume that for $\eta_i$-a.e.\ $y \in Y_i$, the action $G_i \actson (\Xtil_{i,y},\mutil_{i,y})$ is simple. Assume that the groups $G_i$ have trivial center. Let $\gamma_i : \R \times Y_i \to \Lambda_i$ be strict $1$-cocycles.

The actions $\beta(\gamma_i)$ given by \eqref{eq.precise-action-repeated} are conjugate if and only if there exist group isomorphisms $\delta : G_1 \to G_2$, $\rho : \Lambda_1 \to \Lambda_2$ and a $\delta$-conjugacy $\Delta : X_1 \to X_2$ such that, denoting by $\Delta_0 : Y_1 \to Y_2$ the associated isomorphism of flows, the $1$-cocycles $\gamma_2 \circ (\id \times \Delta_0)$ and $\rho \circ \gamma_1$ are cohomologous.
\end{proposition}
\begin{proof}
We start by proving the following claim. Under the assumptions of Proposition \ref{prop.when-induced}, if the group $G$ has trivial center and if $\Sigma \lhd G \times \Lambda$ is a normal subgroup whose action on $X \times \Lambda$ admits a fundamental domain, then $\Sigma \subset \{e\} \times \Lambda$.

Since $\Sigma \actson X \times \Lambda$ admits a fundamental domain, a fortiori, the same holds for the action of $\Sigma$ on $X \times \R \times \Lambda$ given in \eqref{eq.Maharam-of-beta-gamma}. Since $G \actson \Xtil_y$ is simple for a.e.\ $y \in Y$, there is no nontrivial normal subgroup of $G$ whose action on $X \times \R \times \Lambda$ admits a fundamental domain. Hence, $\Sigma \cap (G \times \{e\}) = \{(e,e)\}$. Let now $(g,a) \in \Sigma$ be an arbitrary element. We have to prove that $g = e$. By normality of $\Sigma$, also $(h g h^{-1},a) = (h,e) (g,a) (h,e)^{-1} \in \Sigma$. Since $(g,a) \in \Sigma$, also $(hgh^{-1}g^{-1},e) \in \Sigma$. Since we have proven that $\Sigma$ intersects $G \times \{e\}$ trivially, it follows that $hgh^{-1}g^{-1}=e$ for all $h \in G$. This means that $g$ belongs to the center of $G$, which is assumed to be trivial. This proves the claim.

Now assume that $\Delta_1 : X_1 \times \Lambda_1 \to X_2 \times \Lambda_2$ is a conjugacy w.r.t.\ the group isomorphism $\delta_1 = G_1 \times \Lambda_1 \to G_2 \times \Lambda_2$. Since the action of $\{e\} \times \Lambda_i$ on $X_i \times \Lambda_i$ has $X_i \times \{e\}$ as a fundamental domain, the claim above implies that $\delta_1(\{e\} \times \Lambda_1) = \{e\} \times \Lambda_2$. We define the group isomorphism $\rho : \Lambda_1 \to \Lambda_2$ such that $\delta_1(e,a) = (e,\rho(a))$ for all $a \in \Lambda_1$.

For every $i \in \{1,2\}$, we consider the Maharam extension $G_i \times \Lambda_i \times \R \actson X_i \times \R \times \Lambda_i$ given by \eqref{eq.Maharam-of-beta-gamma}, together with the factor map $\Phi_i : X_i \times \R \times \Lambda_i \to X_i \times \Lambda_i$ defined by \eqref{eq.map-Phi}. Therefore, $\Delta_1$ canonically lifts to a nonsingular isomorphism $\Delta_2 : X_1 \times \R \times \Lambda_1 \to X_2 \times \R \times \Lambda_2$ that is $\R$-equivariant and a $\delta_1$-conjugacy. Since $\Delta_2$ is a $\rho$-conjugacy for the actions of $\Lambda_i$, the map $\Delta_2$ must be of the form $\Delta_2(x,s,b) = (\Delta_3(x,s),\zeta(x,s) \rho(b))$.

Since $\delta_1(e,a) = (e,\rho(a))$ where $\rho : \Lambda_1 \to \Lambda_2$ is an isomorphism, the isomorphism $\delta_1$ must be of the form $\delta_1(g,a) = (\delta(g),\al(g)\rho(a))$, where $\delta : G_1 \to G_2$ is an isomorphism and $\al : G_1 \to \Lambda_2$ is a group homomorphism. Expressing that $\Delta_2((g,e) \cdot (x,s,e)) = (\delta(g),\al(g)) \cdot \Delta_2(x,s,e)$, we find that $\zeta(g \cdot (x,s))^{-1} = \al(g) \zeta(x,s)^{-1}$. By our assumptions, for a.e.\ $y \in Y_1$, the action $G_1 \actson \Xtil_{1,y}$ is not induced. Therefore, $\al(g) = e$ for all $g \in G_1$ and $\zeta(x,s) = \zeta_1(\pi_1(x,s))$ where $\zeta_1 : Y_1 \to \Lambda_2$.

By construction, $\Delta_3 : X_1 \times \R \to X_2 \times \R$ is $\R$-equivariant and a $\delta$-conjugacy. After replacing $\mu_2$ by an equivalent probability measure, we find that $\Delta_3(x,s) = (\Delta(x),s)$, where $\Delta : X_1 \to X_2$ is a measure preserving $\delta$-conjugacy, so that $\om_2(\delta(g),\Delta(x)) = \om_1(g,x)$. The $\delta$-conjugacy $\Delta : X_1 \to X_2$ induces an isomorphism $\Delta_0 : Y_1 \to Y_2$ of associated flows. By construction, $\Delta_0 \circ \pi_1 = \pi_2 \circ \Delta_3$.

Since $\Delta_2$ is $\R$-equivariant, we also find that
$$\zeta_1(t \cdot \pi_1(x,s)) \, \rho(\gamma_1(t,\pi(x,s))) = \gamma_2(t,\pi_2(\Delta(x),s)) \, \zeta_1(\pi(x,s)) \; .$$
This precisely means that the $1$-cocycles $\gamma_2 \circ (\id \times \Delta_0)$ and $\rho \circ \gamma_1$ are cohomologous.
\end{proof}

Proposition \ref{prop.complete-classification} is now an immediate consequence of Corollary \ref{cor.key-OE-superrigid} and Propositions \ref{prop.when-induced} and \ref{prop.when-conjugate}. For completeness, we include a detailed argument.

\begin{proof}[{Proof of Proposition \ref{prop.complete-classification}}]
Point~1 was already proven in Corollary \ref{cor.key-OE-superrigid}. To prove point~2, we start by the easy implication. If we are given all the data mentioned in point~2, we replace $\mu'$ by the equivalent measure $\Delta_*(\mu)$, so that $\om'(\delta(g),\Delta(x)) = \om(g,x)$. We replace $\gamma,\gamma'$ by the cohomologous $1$-cocycles $\gamma_0,\gamma_0'$. Further replacing $\gamma'_0$ by a cohomologous, $\Lambda'_0$-valued $1$-cocycle, we may assume that $\rho \circ \gamma_0 = \gamma'_0 \circ (\id \times \Delta_0)$. By construction, $\Delta \times \rho$ defines a $(\delta \times \rho)$-conjugacy of the actions $\beta(\gamma_0)$ and $\beta(\gamma'_0)$.

Denote by $\si_0, \si'_0$ the restriction of $\si,\si'$ to $G \times \Lambda_0$, $G \times \Lambda'_0$. Since $\al \circ \si_0 = \si'_0 \circ (\delta \times \rho)$, it follows from the previous paragraph that the actions $\beta(\gamma_0,\si_0)$ and $\beta(\gamma'_0,\si'_0)$ are $\al$-conjugate. By construction, the first action is isomorphic with $\beta(\gamma,\si)$ and the second action is isomorphic with $\beta(\gamma',\si')$. We have thus proven that $\beta(\gamma,\si)$ and $\beta(\gamma',\si')$ are conjugate.

Conversely, assume that there is an $\al$-conjugacy between $\beta(\gamma,\si)$ and $\beta(\gamma',\si')$. When an ergodic group action $\Gamma \actson Z$ is induced from both $\Gamma_0 \actson Z_0$ and $\Gamma_1 \actson Z_1$, there exists a $g \in \Gamma$ such that $\Gamma \actson Z$ is induced from the action of $\Gamma_0 \cap g \Gamma_1 g^{-1}$ on $Z_0 \cap g \cdot Z_1$. After composing $\al$ with an inner automorphism and using Proposition \ref{prop.when-induced}, we find subgroups $\Lambda_0 < \Lambda$ and $\Lambda'_0 < \Lambda'$ such that $\al(\si(G \times \Lambda_0)) = \si'(G \times \Lambda'_0)$ and such that $\gamma,\gamma'$ are cohomologous with $\gamma_0,\gamma'_0$ taking values in $\Lambda_0,\Lambda'_0$. We also find the group isomorphism $\delta_1 : G \times \Lambda_0 \to G \times \Lambda'_0$, with $\al \circ \si = \si' \circ \delta_1$, such that the actions $\beta(\gamma_0)$ and $\beta(\gamma'_0)$ are $\delta_1$-conjugate. The conclusion then follows from Proposition \ref{prop.when-conjugate}.
\end{proof}

We finally turn to our main family of group actions, given by \eqref{eq.our-actions}. For $i \in \{1,2\}$, let $G_i \actson (X_i,\mu_i)$ be nonsingular actions and denote by $\om_i : G_i \times X_i \to \R$ the logarithm of the Radon-Nikodym cocycle. Let $\R \actson (Y_i,\eta_i)$ be nonsingular flows. Consider the actions
\begin{equation}\label{eq.our-actions-repeated}
\si_i : G_i \actson X_i \times Y_i : g \cdot (x,y) = (g \cdot x, \om_i(g,x) \cdot y) \; .
\end{equation}
If $\Delta_1 : X_1 \to X_2$ is a conjugacy between the actions $G_1 \actson X_1$ and $G_2 \actson X_2$ and if $\Delta_2 : Y_1 \to Y_2$ is an isomorphism of the flows, writing $\rho(x) = \log(d((\Delta_1)_* \mu_1)/d\mu_2)$, the map
\begin{equation}\label{eq.formula-conjugacy}
X_1 \times Y_1 \to X_2 \times Y_2 : (x,y) \mapsto (\Delta_1(x),-\rho(\Delta_1(x)) \cdot \Delta_2(y))
\end{equation}
is a conjugacy between $\si_1$ and $\si_2$.

We now prove in Proposition \ref{prop.special-when-conjugate} that under the appropriate assumption, the converse also holds. In Lemma \ref{lem.conjugacy}, we explain that this assumption indeed holds for actions of the form $G \actson \R^n/\R^*_+$ when $G < \SL(n,\R)$ is a countable dense subgroup. This then leads to a proof of Corollary \ref{cor.main-example}.

\begin{proposition}\label{prop.special-when-conjugate}
For $i \in \{1,2\}$, let $G_i \actson (X_i,\mu_i)$ be essentially free, ergodic, nonsingular actions of type III$_1$. Assume that the Maharam extensions $G_i \actson (\Xtil_i,\mutil_i)$ admit a unique measure scaling action commuting with the $G_i$-action. Let $\R \actson (Y_i,\eta_i)$ be ergodic flows.

Every conjugacy between the actions $\si_i$ in \eqref{eq.our-actions-repeated} is of the form \eqref{eq.formula-conjugacy} for a conjugacy $\Delta_1$ between $G_1 \actson X_1$ and $G_2 \actson X_2$ and an isomorphism of flows $\Delta_2 : Y_1 \to Y_2$.
\end{proposition}
\begin{proof}
Denote by $\R^2 \actson (Z_i,\zeta_i)$ the unique actions associated with the flows $\R \actson (Y_i,\eta_i)$ given by Proposition \ref{prop.unique-char-adj}, together with the factor maps $\pi_i : Z_i \to Y_i$ satisfying $\pi_i((t,r) \cdot z) = t \cdot \pi_i(z)$. We realize the Maharam extension of $\si_i$ together with its measure scaling action as
\begin{equation}\label{eq.my-maharam}
\sitil_i : G_i \times \R \actson (X_i \times Z_i,\mu_i \times \zeta_i) : (g,r) \cdot (x,z) = (g \cdot x , (\om(g,x),r) \cdot z) \; .
\end{equation}
We have the canonical factor map $\Psi_i : X_i \times Z_i \to X_i \times Y_i : \Psi_i(x,z) = (x,\pi_i(z))$ satisfying $\Psi_i((g,r) \cdot (x,z)) = g \cdot \Psi_i(x,z)$.

Assume that $\delta : G_1 \to G_2$ is a group isomorphism and $\Delta : X_1 \times Y_1 \to X_2 \times Y_2$ is a $\delta$-conjugacy between $\si_1$ and $\si_2$. Denote by $\Deltatil : X_1 \times Z_1 \to X_2 \times Z_2$ the canonical measure preserving lift, which is a $(\delta \times \id)$-conjugacy for the actions $\sitil_i$ in \eqref{eq.my-maharam} and which satisfies $\Psi_2 \circ \Deltatil = \Delta \circ \Psi_1$.

Define the action $\gamma : \R \actson X_i \times Z_i : \gamma_t(x,z) = (x,(t,0) \cdot z)$. Note that $\gamma$ commutes with $\sitil_i$. We claim that $\Deltatil$ is automatically $\gamma$-equivariant.

To prove this claim, we temporarily identify $Z_i = \R \times \Yh_i$ with the action $\R^2 \actson Z_i$ given by $(t,r) \cdot (s,\yh) = (t+r+\beh(r,\yh)+s,r \cdot \yh)$. Under this identification, $\Deltatil : X_1 \times \R \times \Yh_1 \to X_2 \times \R \times \Yh_2$ is a $\delta$-conjugacy for the actions
\begin{equation}\label{eq.hulpactie}
G_i \actson X_i \times \R \times \Yh_i : g \cdot (x,s,\yh) = (g \cdot x, \om_i(g,x) + s, \yh) \; .
\end{equation}
Since $G_i \actson X_i$ is of type III$_1$, the Maharam extension $G_i \actson X_i \times \R$ is ergodic. Therefore, $\Deltatil$ must be of the form $\Deltatil(x,s,\yh) = (\theta_{\yh}(x,s),\Phi(\yh))$, where $\Phi : \Yh_1 \to \Yh_2$ is a nonsingular isomorphism and, for a.e.\ $\yh \in \Yh_1$, the map $\theta_{\yh} : X_1 \times \R \to X_2 \times \R$ is a $\delta$-conjugacy of the Maharam extensions $G_i \actson X_i \times \R$. Since $\Deltatil$ is measure preserving, a.e.\ $\theta_{\yh}$ is measure scaling. We assumed that these Maharam extensions admit a unique commuting measure scaling action. It follows that $\theta_{\yh}$ is equivariant w.r.t.\ translation in the second variable. This means that $\Deltatil$ is equivariant w.r.t.\ translation in the second variable. Thus, the claim is proven.

We thus consider the actions $G_i \times \R^2 \actson X_i \times Z_i : (g,t,r) \cdot (x,z) = (g \cdot x, (\om(g,x)+t,r) \cdot z)$ and we have proven that $\Deltatil$ is a $(\delta \times \id)$-conjugacy between these actions. Since the action of $\R^2$ on $Z_i$ is ergodic, this forces $\Deltatil$ to be of the form $\Deltatil(x,z) = (\Delta_1(x),\Phi_x(z))$, where $\Delta_1 : X_1 \to X_2$ is a $\delta$-conjugacy and, for a.e.\ $x \in X_1$, $\Phi_x : Z_1 \to Z_2$ is an isomorphism between the actions $\R^2 \actson Z_i$.

Define $\rho(x) = \log(d((\Delta_1)_* \mu_1)/d\mu_2)$ and denote $\Phi'_x(z) = (\rho(\Delta_1(x)),0) \cdot \Phi_x(z)$. We still have that $\Phi'_x$ is an isomorphism between the actions $\R^2 \actson Z_i$. W.r.t.\ the measures $\mu_i \times \zeta_i$ on $X_i \times Z_i$, the isomorphism $\Deltatil$ is measure preserving. It then follows that a.e.\ $\Phi'_x$ is measure preserving w.r.t.\ the measures $\zeta_i$ on $Z_i$.

Expressing that $\Deltatil$ is a $\delta$-conjugacy for the actions of $G_i$ and using that $\Phi_x$ is an isomorphism between the actions of $\R^2$, we find that $(\om_1(g,x),0) \cdot \Phi_{g \cdot x}(z) = (\om_2(\delta(g),\Delta_1(x)),0) \cdot \Phi_x(z)$. Since $\om_i$ are the logarithms of the Radon-Nikodym cocycles for $G_i \actson (X_i,\mu_i)$ and since $\Delta_1$ is a $\delta$-conjugacy, we have by definition of $\rho$ that $\om_2(\delta(g),\Delta_1(x)) = \rho(\Delta_1(x)) - \rho(\Delta_1(g \cdot x)) + \om_1(g,x)$. We conclude that $\Phi'_{g \cdot x}(z) = \Phi'_x(z)$. Since the action $G_1 \actson X_1$ is ergodic, we find a measure preserving isomorphism $\Phi : Z_1 \to Z_2$ between the actions $\R^2 \actson Z_i$ such that $\Phi'_x = \Phi$ for a.e.\ $x \in X_1$.

For such an isomorphism $\Phi$, there is a unique isomorphism $\Delta_2 : Y_1 \to Y_2$ for the actions $\R \actson Y_i$ such that $\pi_2 \circ \Phi = \Delta_2 \circ \pi_1$. Define the $\delta$-conjugacy $\Delta_0 : X_1 \times Y_1 \to X_2 \times Y_2$ by \eqref{eq.formula-conjugacy}, i.e.\ $\Delta_0(x,y) = (\Delta_1(x),-\rho(\Delta_1(x)) \cdot \Delta_2(y))$. By construction, $\Psi_2 \circ \Deltatil = \Delta_0 \circ \Psi_1$. It follows that $\Delta = \Delta_0$.
\end{proof}

Recall that a unitary representation $\pi : \cG \to \cU(H)$ of a locally compact group $\cG$ is said to be a $C_0$-representation if for every $\eps > 0$ and $\xi,\eta \in H$, there exists a compact subset $K \subset \cG$ such that $|\langle \pi(g)\xi,\eta\rangle| < \eps$ for all $g \in \cG \setminus K$.

\begin{lemma}\label{lem.C0-rep}
Let $n \geq 1$ be an integer. Write $\cG = \GL(n,\R) \ltimes \R^n$. Then the unitary representation $\pi : \cG \to \cU(L^2(\R^n)) : (\pi(A,a)\xi)(x) = |\det A|^{-1/2} \xi(A^{-1}(x) -a)$ is a $C_0$-representation.

In particular, the action $\cG \actson \R^n : (A,a) \cdot x = A(a+x)$ induces a homeomorphism of $\cG$ onto a closed subgroup of the Polish group of nonsingular automorphisms of $\R^n$ with the Lebesgue measure.
\end{lemma}
\begin{proof}
Assume the contrary. We then find $\xi,\eta \in L^2(\R^n)$, $\eps > 0$ and a sequence $(A_k,a_k) \in \cG$ that tends to infinity in $\cG$ such that $|\langle \pi(A_k,a_k) \xi, \eta\rangle| > \eps$ for all $k$.

We view $\GL(n,\R)$ as a subgroup of $\cG$ and we denote by $D_n < \GL(n,\R)$ the subgroup of diagonal matrices with positive real numbers on the diagonal. Since $\GL(n,\R) = O(n,\R) D_n O(n,\R)$, we can write $(A_k,a_k) = g_k (d_k,b_k) h_k$ with $g_k,h_k \in O(n,\R)$, $d_k \in D_n$ and $b_k \in \R^n$. After passage to a subsequence, we may assume that $g_k$ and $h_k$ converge to $g$, resp.\ $h$. Replacing $\xi$ by $\pi(h) \xi$ and replacing $\eta$ by $\pi(g)^* \eta$, we may then further assume that $|\langle \pi(d_k,b_k) \xi, \eta\rangle| > \eps$ for all $k$.

To reach a contradiction, it thus suffices to prove that the representation $\theta : \R^*_+ \ltimes \R \to \cU(L^2(\R)) : (\theta(d,b)\xi)(x) = d^{-1/2} \xi(d^{-1} x - b)$ is a $C_0$-representation. Denote by $\lambda$ the Lebesgue measure on $\R$. When $N \in \R_+$ and $\cU,\cV \subset [-N,N]$ are Borel sets with indicator functions $1_\cU, 1_\cV \in L^2(\R)$, we have $\theta(d,b) 1_\cU = d^{-1/2} 1_{d(\cU+b)}$, so that $\langle \theta(d,b) 1_\cU, 1_\cV \rangle = d^{-1/2} \, \lambda(d(\cU + b) \cap \cV)$. We conclude that
\begin{multline*}
\bigl\{ (d,b) \in \R^*_+ \times \R \bigm| |\langle \theta(d,b) 1_\cU, 1_\cV \rangle| \geq \eps \bigr\} \\
\subset \bigl\{ (d,b) \in \R^*_+ \times \R \bigm| \eps/2N \leq d^{1/2} \leq 2N/\eps \; , \; |b| \leq (d^{-1} + 1) N \bigr\} \; ,
\end{multline*}
which is compact.
\end{proof}

\begin{lemma}\label{lem.conjugacy}
For $i = 1,2$ and integers $n_i \geq 2$, let $G_i < \SL(n_i,\R)$ be dense subgroups and consider the actions $G_i \actson \R^{n_i}$.

If $\Delta : \R^{n_1} \to \R^{n_2}$ is a $\delta$-conjugacy between these actions, we have $n_1 = n_2$ and there is a unique $A \in \GL(n_1,\R)$ such that $\Delta(x) = A(x)$ for a.e.\ $x \in \R^{n_1}$ and $\delta(g) = AgA^{-1}$ for all $g \in G_1$.

In particular, if $n \geq 2$ and $G < \SL(n,\R)$ is a dense subgroup, then the action $\R \actson \R^n : t \cdot x = e^{-t/n} x$ is the unique measure scaling action that commutes with $G \actson \R^n$.
\end{lemma}

\begin{proof}
Denote by $\Aut(\R^{n_i})$ the Polish group of nonsingular automorphisms of $\R^{n_i}$. By Lemma \ref{lem.C0-rep}, we may view $\SL(n_i,\R)$ as a closed subgroup of $\Aut(\R^{n_i})$. By our assumption, $\Delta G_1 \Delta^{-1} = G_2$. Taking the closure in $\Aut(\R^{n_i})$, we find that $\Delta \SL(n_1,\R) \Delta^{-1} = \SL(n_2,\R)$. This means in particular that $\delta$ extends to a group isomorphism and homeomorphism $\delta : \SL(n_1,\R) \to \SL(n_2,\R)$. It follows that $n_1 = n_2$ and we write $n = n_1 = n_2$.

Write $\cG = \SL(n,\R)$ and denote by $\cH < \cG$ the closed subgroup fixing the first basis vector $e_1$. We view $\Delta$ as a $\delta$-conjugacy for the transitive action $\cG \actson \cG/\cH$. This means that $\delta(\cH) = g_0 \cH g_0^{-1}$ for some $g_0 \in \cG$ and $\Delta(g\cH) = \delta(g)g_0\cH$ for a.e.\ $g \in \cG$. Since $\delta$ is an automorphism of $\SL(n,\R)$, we find $B \in \GL(n,\R)$ such that either $\delta(g) = B g B^{-1}$ for all $g \in \cG$, or $\delta(g) = B (g^{-1})^T B^{-1}$ for all $g \in \cG$. In the second case, $\delta(\cH)$ and $\cH$ are not conjugate. So, we are in the first case. We get that $g_0 = B C$ where $C \in \GL(n,\R)$ normalizes $\cH$. That means that $C e_1 = a e_1$ for some $a \in \R^*$. Translating the formula $\Delta(g\cH) = \delta(g)g_0\cH$ to $\R^n$, we have proven that $\Delta(x) = a B(x)$ for a.e.\ $x \in \R^n$. Writing $A = a B$, the first part of the lemma is proven.

In particular, if $G < \SL(n,\R)$ is a dense subgroup, the only nonsingular automorphisms of $\R^n$ that commute with $G \actson \R^n$ are given by $x \mapsto a x$ for some $a \in \R^*$. This transformation scales the measure by $|a|^n$. The only measure scaling action $\R \actson \R^n$ commuting with $G \actson \R^n$ is thus given by $t \cdot x = e^{-t/n} x$.
\end{proof}

We have now gathered enough material to prove Corollary \ref{cor.main-example}. We also add the computation of the outer automorphism group of these type III$_0$ orbit equivalence relations. For every ergodic flow $\R \actson^\al (Y,\eta)$, we denote by $\Aut_\R(\al)$ the Polish group of all nonsingular automorphisms of $(Y,\eta)$ that commute with the flow $\al$. Note that $\R \subset \Aut_\R(\al)$ by definition.

\begin{corollary}\label{cor.main-example-with-out}
Corollary \ref{cor.main-example} holds.

Moreover, writing $K = \{a/b \mid a,b \in \cA, b \neq 0\}$, the outer automorphism group $\Out(\cR(n,\cA,\al))$ of the orbit equivalence relation $\cR(n,\cA,\al)$ of the action $\be(n,\cA,\al)$ is given by
\begin{equation}\label{eq.outer-aut-general}
\Out(\cR(n,\cA,\al)) \cong \frac{\cN_{\GL(n,K)}(E(n,\cA))}{K^* \cdot E(n,\cA)} \times \{\pm 1\} \times \Aut_\R(\al) \; .
\end{equation}
In particular, when $\cA = \Z[\cS^{-1}]$ for some finite nonempty set of prime numbers $\cS$, we have $E(n,\cA) = \SL(n,\Z[\cS^{-1}])$ and
\begin{equation}\label{eq.outer-aut-specific-ZS}
\Out(\cR(n,\Z[\cS^{-1}],\al)) \cong \bigl(\Z/n\Z\bigr)^{|S|} \times \{\pm 1\} \times \Aut_\R(\al) \; .
\end{equation}
When $\cA = \cO_K$ where $K \subset \R$ is an algebraic number field with $[K:\Q] \geq 2$, we denote by $\Cl(K)$ its ideal class group and consider the subgroup $\Cl_n(K) = \{J \in \Cl(K) \mid J^n = 1\}$. Then $E(n,\cA) = \SL(n,\cO_K)$ and
\begin{equation}\label{eq.outer-aut-specific-OK}
\Out(\cR(n,\cO_K,\al)) \cong \frac{\cO_K^*}{(\cO_K^*)^n} \times \Cl_n(K) \times \{\pm 1\} \times \Aut_\R(\al) \; .
\end{equation}
\end{corollary}
\begin{proof}
By Theorem \ref{thm.rings-cocycle-superrigid}, under the hypotheses of Corollary \ref{cor.main-example}, the actions $E(n,\cA) \actson \R^n / \R^*_+$ satisfy all the assumptions of Corollary \ref{cor.OE-superrigid-prescribed-flow}. So by Corollary \ref{cor.OE-superrigid-prescribed-flow}, the actions $\be(n,\cA,\al)$ are essentially free, ergodic, simple and OE superrigid (v2), with associated flow $\alh$.

In particular, if $\be(n,\cA,\al)$ and $\be(n',\cA',\al')$ are stably orbit equivalent, the actions must be conjugate. By Proposition \ref{prop.special-when-conjugate} and Lemma \ref{lem.conjugacy}, the flows $\al$ and $\al'$ are isomorphic and the actions $E(n,\cA) \actson \R^n/\R^*_+$ and $E(n',\cA') \actson \R^{n'}/\R^*_+$ are conjugate. Then also their Maharam extensions $E(n,\cA) \actson \R^n$ and $E(n',\cA') \actson \R^{n'}$ are conjugate. By Lemma \ref{lem.conjugacy}, $n=n'$ and there is an $A \in \GL(n,\R)$ such that $A E(n,\cA) A^{-1} = E(n,\cA')$.

For every subring $\cA_1 \subset \R$, we denote by $M(n,\cA_1)$ the ring of $n \times n$ matrices with entries in $\cA_1$. We claim that the subring $N(n,\cA_1) \subset M(n,\cA_1)$ generated by $E(n,\cA_1)$ is equal to $M(n,\cA_1)$. For all $i,j \in \{1,\ldots,n\}$ and $a \in \cA_1$, we denote by $E_{ij}(a)$ the matrix that has the entry $a$ in position $ij$ and $0$'s elsewhere. When $i \neq j$ and $a \in \cA_1$, we have that $1 + E_{ij}(a) \in E(n,\cA_1)$ and $1 \in E(n,\cA_1)$. Thus, $E_{ij}(a) \in N(n,\cA_1)$. Since $\SL(n,\Z) = E(n,\Z) \subset E(n,\cA_1)$, given $i \neq j$, the matrix $\si_{ij}$ with entry $1$ in position $ij$, entry $-1$ in position $ji$, and $0$'s elsewhere, belongs to $E(n,\cA_1)$. Thus, also $E_{ii}(a) = \si_{ij} E_{ji}(a) \in N(n,\cA_1)$. This proves the claim.

Since $A E(n,\cA) A^{-1} = E(n,\cA')$, the claim above implies that $A M(n,\cA) A^{-1} = M(n,\cA')$. So, for every $i \in \{1,\ldots,n\}$ and $a \in \cA$, we have that $A E_{ii}(a) A^{-1} \in M(n,\cA')$. Taking the $jj$-entry, it follows that $A_{ji} (A^{-1})_{ij} a \in \cA'$. Summing over $i$, it follows that $a \in \cA'$. So, $\cA \subset \cA'$. By symmetry, also the converse inclusion holds, so that $\cA = \cA'$.

The arguments above apply in particular to the self orbit equivalences of $\be(n,\cA,\al)$ and give us that
$$
\Out(\cR(n,\cA,\al)) \cong \frac{\cN_{\GL(n,\R)}(E(n,\cA))}{\R^*_+ \, E(n,\cA)} \times \Aut_\R(\al) \; .
$$
Write $\cG_+ = \{A \in \GL(n,\R) \mid \det A > 0\}$. Since $n$ is odd, we have that $\GL(n,\R) = \{\pm 1\} \times \cG_+$. It follows that
\begin{equation}\label{eq.intermediate-out}
\Out(\cR(n,\cA,\al)) \cong \frac{\cN_{\cG_+}(E(n,\cA))}{\R^*_+ \, E(n,\cA)} \times \{\pm 1\} \times \Aut_\R(\al) \; .
\end{equation}
Above, we have also proven that every $A \in \cN_{\GL(n,\R)}(E(n,\cA))$ satisfies $A M(n,\cA) A^{-1} = M(n,\cA)$, meaning that $A_{ij} (A^{-1})_{kl} \in \cA$ for all $i,j,k,l$. Denoting by $K = \{a/b \mid a,b \in \cA, b \neq 0\}$ the field of fractions of $\cA$, it follows in particular that $A$ must be a multiple of a matrix with entries in $K$. It follows that
\begin{equation}\label{eq.tussenresultaat}
\frac{\cN_{\cG_+}(E(n,\cA))}{\R^*_+ \, E(n,\cA)} = \frac{\cN_{\cG_+ \cap \GL(n,K)}(E(n,\cA))}{(\R^*_+ \cap K) \, E(n,\cA)} = \frac{\cN_{\GL(n,K)}(E(n,\cA))}{K^* \, E(n,\cA)} \; .
\end{equation}
Combining this with \eqref{eq.intermediate-out}, we have proven \eqref{eq.outer-aut-general}.

By \cite[Theorem B]{OM65}, a matrix $A \in \GL(n,K)$ normalizes $E(n,\cA)$ if and only if $A(\cA^n) = \mathfrak{a} \cA^n$ for an invertible fractional ideal $\mathfrak{a} \subset K$.

When $\cA = \Z[\cS^{-1}]$, every fractional ideal of $\cA$ is principal, so that $\cN_{\GL(n,K)}(E(n,\cA)) = K^* \GL(n,\cA)$. By \cite[Theorem 4.3.9]{HOM89}, we also have that $E(n,\cA) = \SL(n,\cA)$. So, the natural map
$$\GL(n,\cA) \to \frac{\cN_{\GL(n,K)}(E(n,\cA))}{K^* \, E(n,\cA)}$$
is surjective and has kernel $\{A \in \GL(n,\cA) \mid \det A \in (\cA^*)^n\}$. The group of units of $\cA = \Z[\cS^{-1}]$ is the free abelian group generated by $p \in \cS$, so that combining \eqref{eq.intermediate-out} en \eqref{eq.tussenresultaat}, we have proven \eqref{eq.outer-aut-specific-ZS}.

When $K$ is an algebraic number field with $[K:\Q] \geq 2$ and $\cA = \cO_K$, again by \cite[Theorem 4.3.9]{HOM89}, we have that $E(n,\cA) = \SL(n,\cA)$. Denote by $I(K)$ the group of fractional ideals in $K$. As mentioned in \cite[Example 6.5]{OM65}, the following two statements hold. If $A \in \GL(n,K)$, $\mathfrak{a} \in I(K)$ and $A(\cA^n) = \mathfrak{a} \cA^n$, we have $\mathfrak{a}^n = \det(A) \, \cA$. Conversely, if $\mathfrak{a} \in I(K)$, $\al \in K^*$ and $\mathfrak{a}^n = \al \, \cA$, there exists an $A \in \GL(n,K)$ such that $A(\cA^n) = \mathfrak{a} \cA^n$ and $\det A = \al$. For completeness, we provide a more detailed argument. The first implication is contained in \cite[81:7]{OM73}. For the second implication, take $\mathfrak{a} \in I(K)$ and $\al \in K^*$ with $\mathfrak{a}^n = \al \, \cA$. Denote by $e_1,\ldots,e_n$ the standard basis of $K^n$. By \cite[81:5]{OM73}, we find $\mathfrak{b} \in I(K)$ and $B \in \GL(n,K)$ such that $\mathfrak{a} \cA^n = B(\mathfrak{b} e_1 + \cA e_2 + \cdots + \cA e_n)$. By our assumption and \cite[81:7]{OM73}, we have $\al \, \cA = \mathfrak{a}^n = \det(B) \, \mathfrak{b}$. Denote by $D \in GL(n,K)$ the diagonal matrix with $D_{11} = \al \det(B)^{-1}$ and $D_{ii} = 1$ for $i \neq 1$. Since $\mathfrak{b} = D_{11} \, \cA$, writing $A = BD$, we conclude that $A(\cA^n) = \mathfrak{a} \cA^n$. By construction, $\det(A) = \al$.

Define the subgroup $X_n(K)$ of the abelian group $I(K) \times K^*$ by $X_n(K) = \{(\mathfrak{a},\alpha) \in I(K) \times K^* \mid \mathfrak{a}^n = \alpha \cA\}$. Define the subgroup $Y_n(K) \subset X_n(K)$ by $Y_n(K) = \{(\be \cA, \be^n) \mid \be \in K^*\}$. Define $V_n(K) = X_n(K) / Y_n(K)$. To conclude the proof of the corollary, we prove the following two statements.
$$V_n(K) \cong \frac{\cA^*}{(\cA^*)^n} \times \Cl_n(K) \quad \text{and}\quad V_n(K) \cong \frac{\cN_{\GL(n,K)}(\SL(n,\cA))}{K^* \, \SL(n,\cA)} \; .$$
The projection on the first coordinate gives a surjective group homomorphism $V_n(K) \to \Cl_n(K)$ with kernel $\cA^* / (\cA^*)^n$. We prove that this homomorphism is split. Since $\Cl_n(K)$ is a finite abelian group in which the order of every element divides $n$, it suffices to prove that every element of order $k \mid n$ in $\Cl_n(K)$ can be lifted to an element of order $k$ in $V_n(K)$. When $\mathfrak{a} \in I(K)$ and $\mathfrak{a}^k = \beta \cA$ with $\beta \in K^*$, write $n = k m$ and note that $(\mathfrak{a},\beta^m)$ defines such a lift of order $k$. So, $V_n(K) \cong \cA^* / (\cA^*)^n \times \Cl_n(K)$.

For every $(\mathfrak{a},\alpha) \in X_n(K)$, by the discussion above, we can choose $A \in \GL(n,K)$ such that $A(\cA^n) = \mathfrak{a} \cA^n$ and $\det A = \alpha$. This matrix $A$ is uniquely determined up to right multiplication by a matrix in $\SL(n,\cA)$. This realizes a surjective group homomorphism
$$X_n(K) \to \frac{\cN_{\GL(n,K)}(\SL(n,\cA))}{K^* \, \SL(n,\cA)}$$
and the kernel of this homomorphism is by construction equal to $Y_n(K)$. This concludes the proof of the corollary.
\end{proof}

\begin{example}\label{ex.allerlei-voorbeelden}
Let $\cA = \Z[\cS^{-1}]$ where $\cS$ is a finite nonempty set of prime numbers, or let $\cA = \cO_K$ be the ring of integers of an algebraic number field $K \subset \R$ with $[K:\Q] \geq 2$. As mentioned above, we have $E(n,\cA) = \SL(n,\cA)$ for all $n \geq 3$, and a proof can for instance be found in \cite[Theorem 4.3.9]{HOM89}.

By \cite[Corollaries 6.6 and 7.10]{Sus77}, for the same rings $\cA$ as in the previous paragraph and for all integers $0 \leq k \leq s$, we also have that $E(n,\cA_1) = \SL(n,\cA_1)$ when $\cA_1$ is the ring of (Laurent) polynomials over $\cA$, defined by $\cA_1 = \cA[X_1,\ldots,X_k,X_{k+1}^{\pm 1},\ldots,X_s^{\pm 1}]$, and $n \geq 3$. Therefore, whenever $\lambda_1,\ldots,\lambda_s \in \R$ are algebraically independent transcendental numbers, we find that the rings
$$\cA_2 = \Z[\cS^{-1},\lambda_1,\ldots,\lambda_k,\lambda_{k+1}^{\pm 1},\ldots,\lambda_s^{\pm 1}] \quad\text{and}\quad \cA_3 = \Z[\cO_K,\lambda_1,\ldots,\lambda_k,\lambda_{k+1}^{\pm 1},\ldots,\lambda_s^{\pm 1}]$$
also satisfy $E(n,\cA_i) = \SL(n,\cA_i)$ for all $n \geq 3$ and $i \in \{2,3\}$. The associated group actions $\be(n,\cA_i,\al)$ are OE superrigid (v2). Moreover the outer automorphism groups $\Out(\be(n,\cA_i,\al))$ are still given by formulas similar to \eqref{eq.outer-aut-specific-ZS} and \eqref{eq.outer-aut-specific-OK}~:
\begin{equation*}\begin{split}
& \Out(\cR(n,\cA_2,\al)) \cong \bigl(\Z/n\Z\bigr)^{|S|+s-k} \times \{\pm 1\} \times \Aut_\R(\al) \; ,\\
& \Out(\cR(n,\cA_3,\al)) \cong \bigl(\Z/n\Z\bigr)^{s-k} \times \frac{\cO_K^*}{(\cO_K^*)^n} \times \Cl_n(K) \times \{\pm 1\} \times \Aut_\R(\al) \; .
\end{split}\end{equation*}
The reason for this is that by \cite[Corollary 5.6]{AA81}, the invertible fractional ideals in the (Laurent) polynomial rings $\cA_1$ above are all the product of a principal ideal and a fractional ideal in $\cA$. The only difference compared to \eqref{eq.outer-aut-specific-ZS} and \eqref{eq.outer-aut-specific-OK} thus comes from the group of units in $\cA_i$.
\end{example}

\begin{remark}\label{rem.a-lot}
We remark that Corollary \ref{cor.main-example} provides a family of OE superrigid actions, with a prescribed associated flow, that is large and complex in a descriptive set theoretic sense of the word. Fix a prime number $p$. For every finite subset $\cF \in \R$, we consider the ring $\Z[p^{-1},\cF]$, which satisfies the assumptions of Corollary \ref{cor.main-example}, so that the actions $\be(n,\Z[p^{-1},\cF],\al)$ are OE superrigid (v2).

By Corollary \ref{cor.main-example}, to decide when two such actions are stably orbit equivalent, we have to decide if $\Z[p^{-1},\cF] = \Z[p^{-1},\cF']$. This defines a complicated equivalence relation $\cR$ on the Borel space of finite subsets of $\R$. When $\lambda,\lambda' \in \R$ are transcendental, we get that $\Z[p^{-1},\lambda] = \Z[p^{-1},\lambda']$ if and only if there exist $a \in \Z[p^{-1}]^*$ and $b \in \Z[p^{-1}]$ with $\lambda = a \lambda' + b$. So at least, the equivalence relation $\cR$ is not smooth.
\end{remark}

\begin{remark}\label{rem.W-star-superrigid}
As mentioned in the introduction, we can combine the construction of \cite[Proposition D]{Vae13} with Corollary \ref{cor.OE-superrigid-prescribed-flow} to give ad hoc examples of nonsingular actions $G \actson (X,\mu)$ that are W$^*$-superrigid (v2) and that have any prescribed associated flow.

As in \cite[Proposition D]{Vae13}, denote by $\Sigma < \SL(5,\Z)$ the subgroup of matrices $A$ satisfying $A(e_i) = e_i$ for $i=1,2$. Define $G$ as the amalgamated free product $G = \SL(5,\Z) \ast_\Sigma (\Sigma \times \Z)$, with canonical homomorphism $\pi : G \to \SL(5,\Z)$. Consider the probability measure preserving Bernoulli action $G \actson ([0,1],\lambda)^G$. Given any ergodic flow $\R \actson^\al (Y,\eta)$, consider the action
\begin{multline*}
\quad\quad \beta_\al : G \actson (\R^5 \times Y)/\R \times [0,1]^G : g \cdot (\overline{(x,y)},z) = (\overline{(\pi(g)x,y)},g \cdot z)\\
\text{where}\quad \R \actson \R^5 \times Y : t \cdot (x,y) = (e^{t/5}x, t \cdot y) \; .\quad\quad
\end{multline*}
Then, $\beta_\al$ is essentially free, nonsingular, ergodic, simple and W$^*$-superrigid (v2), with associated flow $\alh$.

To prove this result, one uses \cite[Theorem 8.1 and Proposition D]{Vae13} to show that the crossed product factor associated with $\beta_\al$ has a unique Cartan subalgebra, up to unitary conjugacy. It then suffices to prove that $\beta_\al$ is OE superrigid (v2). Since $G$ is finitely generated and has trivial center, by Corollary \ref{cor.OE-superrigid-prescribed-flow}, it suffices to prove that the action $G \actson \R^5 \times [0,1]^G : g \cdot (x,z) = (\pi(g) x, g \cdot z)$ is simple and cocycle superrigid with countable targets. Simplicity is easy to check and cocycle superrigidity was proven in \cite[Proposition D]{Vae13}.
\end{remark}


\begin{thebibliography}{Abcd12}\setlength{\itemsep}{-1mm} \setlength{\parsep}{0mm} \small
\bibitem[AA81]{AA81} D.D. Anderson and D.F. Anderson, Divisorial ideals and invertible ideals in a graded integral domain. {\it J. Algebra} {\bf 76} (1982), 549-569.

\bibitem[BHV08]{BHV08} B. Bekka, P. de la Harpe and A. Valette, Kazhdan's property (T). {\it New Mathematical Monographs} {\bf 11}, Cambridge University Press, Cambridge, 2008.

\bibitem[DIP19]{DIP19} D. Drimbe, A. Ioana and J. Peterson, Cocycle superrigidity for profinite actions of irreducible lattices. To appear in {\it Groups Geom. Dyn.} \href{https://arxiv.org/abs/1910.08642}{arXiv:1910.08642}

\bibitem[DV21]{DV21} D. Drimbe and S. Vaes, Superrigidity for dense subgroups of Lie groups and their actions on homogeneous spaces. To appear in {\it Mathematische Annalen.} \href{https://arxiv.org/abs/2107.06159}{arXiv:2107.06159}

\bibitem[FMW04]{FMW04} D. Fisher, D.W. Morris and K. Whyte, Nonergodic actions, cocycles and superrigidity. {\it New York J. Math.} {\bf 10} (2004), 249-269.

\bibitem[Fur98]{Fur98} A. Furman, Orbit equivalence rigidity. {\it Ann. of Math.} {\bf 150} (1999), 1083-1108.

\bibitem[HOM89]{HOM89} A.J. Hahn and O.T. O'Meara, The classical groups and $K$-theory. {\it Grundlehren der mathematischen Wissenschaften} {\bf 291}. Springer-Verlag, Berlin, 1989.

\bibitem[Ioa08]{Ioa08} A. Ioana, Cocycle superrigidity for profinite actions of property (T) groups. {\it Duke Math. J.} {\bf 157} (2011), 337-367.

\bibitem[Ioa14]{Ioa14} A. Ioana, Strong ergodicity, property (T), and orbit equivalence rigidity for translation actions. {\it J. Reine Angew. Math.} {\bf 733} (2017), 203-250.

\bibitem[JS85]{JS85} V.F.R. Jones and K. Schmidt, Asymptotically invariant sequences and approximate finiteness. {\it Amer. J. Math.} {\bf 109} (1987), 91-114.

\bibitem[OM65]{OM65} O.T. O'Meara, The automorphisms of the linear groups over any integral domain. {\it J. Reine Angew. Math.} {\bf 223} (1966), 56-100.

\bibitem[OM73]{OM73} O.T. O'Meara, Introduction to quadratic forms. Third corrected printing. {\it Die Grundlehren der mathematischen Wissenschaften} {\bf 117}. Springer-Verlag, Berlin, Heidelberg, New York, 1973.

\bibitem[PR94]{PR94} V. Platonov and A. Rapinchuk, Algebraic groups and number theory. {\it Pure and Applied Mathematics} {\bf 139}. Academic Press, Boston, 1994.

\bibitem[Pop05]{Pop05} S. Popa, Cocycle and orbit equivalence superrigidity for malleable actions of $w$-rigid groups. {\it Invent. Math.} {\bf 170} (2007), 243-295.

\bibitem[Pop06]{Pop06} S. Popa, On the superrigidity of malleable actions with spectral gap. {\it J. Amer. Math. Soc.} {\bf 21} (2008), 981-1000.

\bibitem[PV08]{PV08} S. Popa and S. Vaes, Cocycle and orbit superrigidity for lattices in $\SL(n,\R)$ acting on homogeneous spaces. In {\it Geometry, rigidity, and group actions}, Chicago Lectures in Math., Univ. Chicago Press, Chicago, 2011, pp.\ 419-451.

\bibitem[PV09]{PV09} S. Popa and S. Vaes, Group measure space decomposition of II$_1$ factors and W$^*$-superrigidity. {\it Invent. Math.} {\bf 182} (2010), 371-417.

\bibitem[Sch79]{Sch79} K. Schmidt, Asymptotically invariant sequences and an action of $\SL(2,\Z)$ on the $2$-sphere. {\it Israel J. Math.} {\bf 37} (1980), 193-208.

\bibitem[Sus77]{Sus77} A.A. Suslin, The structure of the special linear group over rings of polynomials. {\it Izv. Akad. Nauk SSSR Ser. Mat.} {\bf 41} (1977), 221-238.


\bibitem[Tak03]{Tak03} M. Takesaki, Theory of operator algebras, II. {\it Encyclopaedia of Mathematical Sciences} {\bf 125}, Springer-Verlag, Berlin, 2003.

\bibitem[Vae13]{Vae13} S. Vaes, Normalizers inside amalgamated free product von Neumann algebras. {\it Publ. Res. Inst. Math. Sci.} {\bf 50} (2014), 695-721.

\bibitem[Zim84]{Zim84} R.J. Zimmer, Ergodic theory and semisimple groups. {\it Monographs in Mathematics} {\bf 81}, Birkh\"{a}user Verlag, Basel, 1984.
\end{thebibliography}
\end{document}